\documentclass[11pt,reqno]{amsart}
\usepackage{amsmath,amssymb,amscd,mathrsfs}

\def\st{Suppose that }
\newtheorem{propo}{Proposition}[section]

\newtheorem{lemma}[propo]{Lemma}
\newtheorem{corol}[propo]{Corollary}

\newtheorem{theo}[propo]{Theorem}

\newcommand{\ld}{,\ldots ,}

\newcommand{\ra}{ \rightarrow }

\newcommand{\lan}{ \langle }
\newcommand{\ran}{ \rangle }

\newcommand{\soc}{\mathop{\rm Soc}\nolimits}
\newcommand{\Soc}{\mathop{\rm Soc}\nolimits}

\newcommand{\ch}{\mathop{\rm char}\nolimits}
\newcommand{\diag}{\mathop{\rm diag}\nolimits}

\newcommand{\Id}{\mathop{\rm Id}\nolimits}
\newcommand{\Irr}{\mathop{\rm Irr}\nolimits}
\newcommand{\irr}{\mathop{\rm Irr}\nolimits}

\newcommand{\CC}{\mathop{\mathbb C}\nolimits}

\newcommand{\ZZ}{\mathop{\mathbb Z}\nolimits}

\newcommand{\al}{\alpha}

\newcommand{\ep}{\varepsilon}
\newcommand{\lam}{\lambda }
\newcommand{\om}{\omega }

\newcommand{\TC}{{\mathbf T}}

\newcommand{\up}{^{-1}}

\newcommand{\si}{\sigma }

\def\d12{{_{12}}}

\def\acf{{algebraically closed field }}

\def\au{{automorphism }}

\def\ccc{{constituent }}

\def\f{{following }}

\def\ho{{homomorphism }}

\def\hot{{homomorphism. }}

\def\ii{{if and only if }}

\def\ir{{irreducible }}

\def\irt{{irreducible. }}
\def\irr{{irreducible representation }}

\def\itf{{It follows that }}

\def\mult{{multiplicity }}

\def\po{{polynomial }}

\def\rep{{representation }}
\def\reps{{representations }}

\def\syl{{Sylow $p$-subgroup }}

\def\SL{{\rm SL}}

 \newcommand{\GC}{{\mathbf G}}

\newcommand{\med}{\medskip}

\begin{document}

\title[]{Low dimensional projective indecomposable modules for
Chevalley groups in defining characteristic}

\author{A. E. Zalesski}

\address{Dipartimento di Matematica e Applicazioni,
Universit\`a degli Studi di Milano-Bicocca, via R. Cozzi 53,
20126, Milano, Italy}


\email{alexandre.zalesski@gmail.com}

\maketitle

\begin{abstract}The paper studies  lower bounds for the dimensions
of projective indecomposable modules for Chevalley groups $G$ in
defining characteristic $p$. The main result extending earlier one
by Malle and Weigel (2008) determines the modules in question of
dimension equal to the order of a Sylow $p$-subgroup of $G$. We
also substantially generalize a result by Ballard (1978) on lower
bounds for the dimensions of projective indecomposable modules and
find lower bounds in some cases where Ballard's bounds are
vacuous.
\end{abstract}

\section{Introduction}
\def\ch{characteristic }
Let $G$ be a finite group  and $p$ a prime number. Let $F$ be a
field of \ch $p$ dividing the order $|G|$ of $G$. We always assume
that $F$ contains a primitive $m$-root of unity, where
$m=|G|/|G|_p$ and $|G|_p$ is the order of a \syl $U$ of $G$. The
group algebra $FG$ of $G$ over $F$ viewed as a left $FG$-module is
called the regular module and  indecomposable direct
summands  are called {\it principal indecomposable $FG$-modules},
customarily abbreviated as PIM's. These are classical objects of
study in the modular \rep theory of finite groups \cite{BN}. One
of the main open problems is to determine their dimensions or at
least provide satisfactory information on the dimensions.

This paper studies this problem for finite Chevalley groups in
defining \ch $p$. For certain groups of small rank the PIM
dimensions have been explicitly computed, but not much is known in
general. A rather detailed survey of the current state of the
problem and the methods used to attack it is done by J. Humphreys
\cite{Hub}.

The results of this paper concentrate mainly on obtaining lower
bounds for the PIM dimensions. The absolute lower bound for a PIM
dimension for any finite group $G$ is $|G|_p$, and this bound is
attained for every Chevalley group. Indeed, every such group has
an \ir $FG$-module of dimension $|G|_p$, hence of $p$-defect 0,
known as the Steinberg module. This is unique if $G$ is
quasi-simple, and this is a PIM. Our use of the term `lower bound'
assumes that we exclude \ir modules of $p$-defect 0. 

The earliest result on  lower bounds for PIM dimensions for
Chevalley groups is due to Ballard \cite[Corollary 5.4]{Ba}, who
considered only non-twisted groups. His result is stated in the
same shape in \cite[\S 9.7]{Hub}. We show in Section 4 how to
extend Ballard's result for arbitrary twisted groups, and we also
obtain a version for it for a parabolic subgroup in place of a
Borel subgroup in the original Ballard's statement.

Ballard's lower bound is not available for every PIM, for
instance, it is useless for any PIM for any non-twisted Chevalley
group over the field of order 2. Therefore,  it is essential to
make it clear when Ballard's type bound is applicable, and this is
desirable to be made in terms of the standard parameterization of
PIM's, specifically, in terms of their socles. Recall that every
PIM (for any finite group $G$) has an \ir socle, which  determines
the PIM. This establishes a bijection between PIM's and the \ir
\reps of $G$, which we refer to here as the standard
parametrization. The PIM corresponding in this way to the trivial
$FG$-module $1_G$ is called 1-PIM in \cite{MW}.

No further result was known over almost 30 years until Malle and
Weigel \cite{MW} determined the 1-PIM's of dimension $|G|_p$. They
did so for all simple groups $G$ and for all primes dividing
$|G|$. 
For Chevalley groups in characteristic $p$ they suggested a method
called the parabolic descent in \cite{MW}. This allows to bound
from below the 1-PIM dimension in terms of   Levi subgroups of
parabolic subgroups of $G$. The method in its original shape does
not work for other PIM's. In this paper we develop the method
further to a level when it can be used for arbitrary PIM's (in
characteristic $p$). This allows to extend the above mentioned
result by Malle and Weigel \cite{MW} as follows:

\begin{theo}\label{wt1}
Let $G$ be a quasi-simple Chevalley group in defining
characteristic $p$, and let $\Phi$ be a $p$-modular PIM of
non-zero defect. Then $\dim \Phi>|G|_p$, unless $\Phi$ is a
$1$-PIM and $G/Z(G)\cong PSL(2,p)$ or ${}^2G_2(3)\cong {\rm
Aut}\,SL(2,8)$.
\end{theo}

The parabolic descent reduces the proof to groups of BN-pair rank
at most 2, and for most of them Theorem \ref{wt1} is already known
to be true. However, for the groups $G=SU(4,p), {}^3D_4(p)$ and
${}^2G_2(3^k)$  the PIM dimensions are not known, although the
character tables are known. These are  not sufficient to rule out
the above three groups, and the parabolic descent method is only
partially helpful. In order to prove Theorem \ref{wt1} for these
groups we also make use of the following simple observation (Lemma
\ref{b1}): $\dim\Phi= |G|_p\cdot (\chi,1_U^G)=|G|_p\cdot
(\chi,\Gamma)$, where $\chi$ is the character of $\Phi$, $\Gamma$
is a Gelfand-Graev character and $U$ is the \syl of $G$. (Here
$1_U$ is the trivial \rep of $U$ and $1_U^G$ is the induced
representation.)

I conjecture that for classical groups $G$ of rank $n$ the 
 dimension of a PIM other than the Steinberg one is at least
 $(n-1)\cdot |G|_p$. Some
progress is achieved in this paper by using a new idea based on
the analysis of common \ir constituents of the ordinary character
of a PIM and the induced module $1_U^G$, where $U$ is a \syl of
$G$. Let $W$ be the Weyl group of $G$ viewed as a group with
$BN$-pair. In favorable circumstances, in particular, for groups
$SL(n,q)$, $n>4$, $E_6(q)$, $E_7(q)$, $E_8(q)$  
the PIM dimension is shown to be at least $d\cdot |G|_p$, where
$d$ is the minimum dimension of a non-linear \irr of $W$ (Theorem
\ref{th1}). If $G=SL(n+1,q)$, $n>3$ then the rank of $G$ is $n$,
$W\cong S_{n+1}$, the symmetric group, and $d=n$. So in this case
the conjecture is confirmed. Note that if $q=2$ then there is a
PIM of dimension $n\cdot |G|_p$; if $q>2$ and $n>2$ then there is
a PIM of dimension $(n+1)\cdot |G|_p$ \cite{Tsu}.

\med
{\bf Notation}. ${\mathbb Q}$, ${\mathbb C}$ are the rational and
complex number fields, respectively, and ${\mathbb Z}$ is the ring
of integers. $F_q$ is the finite field of $q$ elements, and $F$ an
\acf of characteristic $p>0$.

If $G$ is a finite group, then $Z(G)$ is the center and $|G|$ is
the order of $|G|$. If $p$ is a prime then $|G|_p$ is the $p$-part
of $|G|$ and also the order of every \syl of $G$. A $p'$-group
means a finite group with no element of order $p$. If $g\in G$
then $|g|$ is the order of $g$. 
Notation for
classical groups are standard.

All \reps and modules are over $F$ or over $\CC$ (unless otherwise
stated). Sometimes we deal with $RG$-modules, where $1\in R\subset
\CC$ is a subring. In this case all $RG$-modules are assumed to be
free as $R$-modules. If the ground field is clear from the
context, we take liberty to use the term `$G$-module'.  All
modules are assumed to be finitely generated. Representations of
$G$ in characteristic $p$ are called $p$-{\it modular}, and those
over the complex numbers are called {\it ordinary}. The regular
\rep of $G$ is denoted by $\rho_G^{reg}$,  and the trivial
one-dimensional \rep is denoted by $1_G$. We also use $1_G$ to
denote the trivial one-dimensional module and its (Brauer)
character. If $M$ is an $FG$-module, then $\soc M$ means the socle
of $M$, the sum of all \ir submodules.

The set of \ir characters of $G$ is denoted by $\Irr G$, and
 $\ZZ \Irr G$ is the $\ZZ$-span of $\Irr G$; elements of  $\ZZ \Irr G$
 are called {\it generalized  characters}.

A projective indecomposable $FG$-module is called a PIM and
usually denoted by $\Phi$.  Every PIM is determined by its socle.
The PIM whose socle is $1_G$ is called here 1-PIM, and denoted by
$\Phi_1$. More notation concerning projective modules is
introduced in Section 3. We set $c_\Phi=\dim \Phi /|G|_p$. 

If $\chi$ is a character vanishing at all non-identity
$p$-elements then we write $c_\chi =\chi(1)/|G_p|$; this is an
integer.

Let $M$ be an $FG$- or $RG$-module. We set $C_M(G)=\{m\in M:gm=m$
for all $g\in G \}$. Thus,  $C_M(N)$ is the set of $G$-invariants
(or fixed points) on $M$. 

Let $H$ be a subgroup of $G$ and $N$ a  normal subgroup of $H$.
Then $C_M(N)$ is an $H$-module, and when it is viewed as an
$H/N$-module, it is denoted  by $\overline{M}_{H/N}$ (or
$\overline{M}$), and called the {\it Harish-Chandra restriction of
$M$ to} $H/N$. Conversely, given an $H/N$-module $D$, one can view
it as an $H$-module with trivial action of $N$. Then the induced
$G$-module $D^G$ (when $D$ is viewed as an $H$-module) is denoted
by $D^{\# G}$ and  called {\it Harish-Chandra induced from} $D$.
For details see  \cite[p. 667-668, \S 70A]{CR2}, where these
operations are called generalized restriction and induction. This
corresponds to similar operations on characters and Brauer
characters. The Harish-Chandra restriction and induction extends
by linearity to the class functions on $G$. So if $\chi$ is a
class function on $G$ then $\overline{\chi}_{H/N}$ is the
corresponding Harish-Chandra restriction of $\chi$ to  $H/N$, and
if $\lam$ is a class function on $H/N$ then $\lam^{\# G}$ denotes
the Harish-Chandra induced class function on $G$. (For the
ordinary induction we use notation $\lam^G$.) Let $\eta:G\ra
{\mathbb C}$ be a class function on $G$. The formula $(\lam^{\#
G},\eta)=(\lam , \overline{\eta}_{H/N})$ is an easy consequence of
the Frobenius reciprocity and called the {\it Harish-Chandra
reciprocity}. (This is formula 70.1(iii) in \cite[p.668]{CR2}.)
Note that $\overline{\eta}_{H/N}$ can be obtained as the
truncation of $\eta$. This is defined by $\frac{1}{|N|}\sum_{n\in
N}\eta(hn)$
for $h\in H$ viewed as a function on $H/N$. 

Let $\GC$  be a connected reductive algebraic group. An algebraic
group \ho $\GC \ra\GC$ is called {\it Frobenius} if its fixed
point subgroup is finite. We denote Frobenius endomorphisms by
$Fr$. So $\GC^{Fr}=\{g\in \GC: Fr(g)=g \}$ is a finite group,
called here a finite reductive group. If $\GC$ is simple and
simply connected, we refer to $G=\GC^{Fr}$ as a {\it Chevalley
group}. More notation concerning algebraic groups will be
introduced in Section 5.
If $G$ is a finite reductive group or an algebraic group, $p$ is
reserved for the defining characteristic of $G$.

\section{The methods and main results}

In order to guide a reader through the paper, we comment here the
machinery  used in the proofs.

\subsection{Parabolic descent}

A well known standard fact on PIM's for any finite group $G$ is
that the mapping $\Phi\ra \Soc \Phi$  yields a bijection
PIM$\,{}_G \ra\Irr G$ between the set of PIM's and the set of
equivalence classes of \ir $FG$-modules. In addition, every
projective module is a sum of PIM's. Therefore,
 every projective module is determined by its socle,
 and a projective module is a PIM \ii its socle is irreducible.

\med
Let $G$ be a Chevalley group, so $G=\GC^{Fr}$, where $\GC$ is a
simply connected simple algebraic group. The \f result is the well
known Smith-Dipper theorem (see \cite{Sm}, \cite{Di}, \cite{Cab},
\cite[2.8.11]{GLS}):

\begin{lemma}\label{sdt}
Let $G$ be a group with a split BN-pair in characteristic $p$, $P$
a parabolic subgroup with Levi subgroup $L$. Let $V$ be an \ir
$FG$-module. Then $\overline{V}_L=C_V(O_p(P))|_L$ is an \ir
$L$-module.
\end{lemma}

In other words, for every parabolic subgroup $P$ of $G$ there is a
 mapping $\sigma_{G,P}:\Irr G\ra 
 \Irr L$ defined via the Harish-Chandra restriction $V\ra \overline{V}_L
$. Note that $C_V(O_p(P))$ coincides with $\soc V|_P$ as $p={\rm char}\, F$. The
 mapping $\sigma_{G,P}$ is surjective (see Lemma \ref{sm9}). 
This allows one to define a surjective mapping
$\pi_{G,P}:\,$PIM${}_G\ra \,$PIM$\,_L $ as the composition of the
mappings

$$\Phi\ra \soc \Phi\ra \overline{(\soc \Phi)}_L\ra \Psi,$$
\smallskip
where $\Psi$ is the PIM for $L$ with socle $\overline{(\soc
\Phi)}_L$. This is well defined in view of Lemma \ref{sdt}. We
call this mapping the {\it parabolic descent from $G$ to $L$}. 
(We borrow the term from \cite{MW} but the meaning of the term is
not the same as in \cite{MW}.) One observes that $\pi_{G,L}(\Phi)$
is a direct summand of $\overline{ \Phi}_L=C_{\Phi}(O_p(P))|_L$
(but the equality  rarely holds). This implies that $c_\Phi\geq
c_{\Psi}$, where $\Psi=\pi_{G,P}(\Phi)$, see Lemma \ref{mp1}. In
the special case where $\Phi=\Phi_1$ is the 1-PIM, this fact has
been exploited in \cite{MW}.    An attempt to extend it to other
PIMs   meets an obstacle, specifically, in order the method could
work one needs at least to guarantee that $\Psi$ is not a defect zero \ir $FL$-module.
We show how to manage with this difficulty in Section 5.

A weak point of the parabolic descent is that it only allows to
bound (from below) $\dim \Phi$ in terms of $\dim\Psi$, and can not
 be used for showing that $c_\Phi$ grows as the rank of $G$ tends
to infinity. Nonetheless this is useful for classifying PIM's of
dimension $|G|_p$, which is one of our tasks below.

\med
 A formal analog of Lemma \ref{sdt} for a projective module would be a
claim that if $\Phi$ is a PIM for $G$ then $\overline{\Phi}_L$ is a
PIM for $L$. However, this is not true. This is evident from
Propositions \ref{r2} and \ref{r55}.

\subsection{Ballard's bound revised}\label{sub1}

Let $B,N$ be the subgroups of $G$ defining the $BN$-pair structure
on $G$ (see \cite[\S 69.1]{CR2}). The group $T:=N\cap B$ is normal
in $N$. Set $W=N/T$. Then $B$ is a Borel subgroup of $G$,
$U=O_p(B)$ is a \syl of $G$ and $T=N\cap B$ is a maximal torus of
$B$. Let $\Phi$ be a PIM for $G$ with socle $M$. Ballard
\cite[Corollary 5.4]{Ba} proves that $\dim \Phi\geq |W\beta|\cdot |G|_p$, equivalently,
$c_\Phi\geq |W\beta|$ in our notation, where $ |W\beta|$ is  the
size of the $W$-orbit of the (Brauer) character $\beta$ of $T$
afforded by the restriction of $\soc (M|_B)$ to $T$. He assumes
$G$ to be non-twisted. In Section 4 Ballard's result is
generalized to all twisted groups as follows. 
The conjugacy action of $N$ on $T$ induces an action on $\Irr T $,
and for $\beta\in\Irr T$ let $|N\beta|$ denote the size of the
$N$-orbit of $\beta$. Let $\Phi$ be a PIM for $G$ with socle $V$
and let $\beta $  be the Brauer character of the $FG$-module
$C_V(U)|_T$; it is well known that $\beta\in\Irr T$. Then our
version of Ballard's result states that $c_\Phi\geq |N\beta|$
(Proposition \ref{r2}).

The bound $|N\beta|$ is vacuous if $ |N\beta|=1$. For instance,
there are PIM's for which $\beta=1_{T}$, so the bound is vacuous
for such PIM's. In addition, in many cases $ |N\beta|$ is small,
and the bound is not sharp. (This happens for instance if
$|N/T|=2$ but $G$ is not $SL(2,q)$.)

\med
We find out that the nature of Ballard's result is not specific
for PIM's. We prove a similar result for arbitrary characters
vanishing at all non-trivial $p$-elements, see Proposition
\ref{r2b}. This implies the result for PIM's as their characters
have this property.

\med
 The parabolic descent can be combined with our
interpretation of the Ballard bound as follows. Let $P$ be a
standard parabolic subgroup of $G$ \cite[65.15]{CR2}, and let $L$
be the standard Levi subgroup of $P$ \cite[69.14]{CR2}. Let $\Psi$
be the parabolic descent of $\Phi$. Let $N_L$
be the normalizer of $L$ in $N$, so $N_L$ acts on $L$ via
conjugation. For $n\in N_L$ let $\Psi^n$ denote the twist of
$\Psi$ by $n$. (For any $FL$-module $X$ one defines $X^n$ to be
$X$ with the twisted action of $L$, that is, $l^n(x):=nln\up\cdot
x$, where $l\in L$, $x\in X$, $n\in N_L $.) Denote by $|N\Psi|$
the size of the $N$-orbit of $\Psi$. Then our
generalization of Ballard's theorem asserts that $c_\phi\geq
|N\Psi|\cdot c_\Psi$ contains every PIM $ \Psi^n$ $(n\in N_L)$ (Proposition
\ref{r55}).

\subsection{Harish-Chandra induction} Let $G\in \{SL(n,q),~ n>4,~
E_6(q),~E_7(q),~E_8(q)\}$, and let $r$ be the rank of $G$. Let
$\chi$ be the character of a PIM $\Phi\neq St$. We show (Section
6) that the Harish-Chandra theory together with the main result of
Malle-Weigel \cite{MW} yields a lower bound $c_{\Phi}\geq r$. Here
is the idea of the proof. In notation of Section \ref{sub1}, by
analysis of the action of $N$ on $T$ (Proposition \ref{r4})
and using Ballard's bound, we deduce that $c_{\Phi}\geq
|N\beta|>r$, whenever $ \beta\neq 1_{T}$. If $1_{T}= \beta$
then Ballard's bound is vacuous. In this case we first show that
$c_\Phi=(\chi,1_B^G)$,  see Proposition \ref{r01a}. Let
$\lam\in\Irr G $  be a common constituent of $\chi$ and $1_B^G$.
The Harish-Chandra
 theory tells us that 
$(\lam, 1_B^G)\geq r$  
 for the
above groups, unless  $\lam\in \{St, 1_G\}$. As $St$ is not a
constituent of $\chi$,  $c_\Phi< r$ implies $\lam=1_G$.  So
$c_\Phi=(\chi,1_B^G)=(\chi,1_G)$. 
By general modular representation theory, $(\chi,1_G)\neq 0$
implies $\Phi=\Phi_1$ and $ (\chi,1_G)=1$, and hence $c_\Phi=1$.
The groups $G$ for which $c_{\Phi}=1$ have been determined in
\cite{MW}.

We expect that this reasoning can be improved to obtain a lower
bound for all classical groups, however, this requires much deeper
analysis.

\section{Preliminaries}

 Let $G$ be a finite group of order $|G|$
and $p$ a prime number. Let  $\ep$ be a primitive $|G|$-root of
unity. Any ordinary representation is equivalent to a
representation $\phi$ over ${\mathbb Q}(\ep)$, and moreover, over
a maximal subring $R $ of ${\mathbb Q}(\ep)$ not containing $1/p$.
In addition, $R$ has a unique maximal ideal $I$ such that $F=R /I$
is a finite field of characteristic $p$. Note that $F$ contains a
primitive $m$-root of unity, where $m=|G|/|G|_p$. For uniformity
one can similarly define $R$ in the algebraic closure of ${\mathbb
Q}$, and then fix this $R $ to deal with all finite groups. The
mapping $R \ra F$ yields also a surjective \ho of the group of
roots of unity in $R $ to the group of roots of unity in $F$, used
to define Brauer characters.

Every ordinary representation is equivalent to a representation
over $R$.  So if $\phi (G)\subset GL(n,R)$ for some $n$, then the
natural projection $GL(n,R)\ra GL(n,F)$ yields  a $p$-modular
representation $\overline{\phi}:G\ra GL(n,F)$ called the {\it
reduction of $\phi$ modulo} $p$. It is well known that the
composition factors of  $\overline{\phi}$ remain \ir under any
field extension of $F$. This can be translated to the language of
$RG$- and $FG$-modules, however, it requires to consider only $R
G$-modules that are free as $R$-modules. So {\it all $RG$-modules
below are assumed to be free as $R$-modules}.

If $G$ is a $p'$-group then, by Dickson's theorem, the reduction
yields a bijection between the isomorphism classes of $ RG$- and
$FG$-modules which makes identical the $p$-modular \rep theory
with the ordinary \rep theory of $G$. If $p$ divides $|G|$, this
is not true anymore, however, there is a rather sophisticated
replacement: the reduction modulo $p$ yields a bijection between
the isomorphism classes of projective $RG$- and projective
$FG$-modules (Swan' theorem, see \cite[Theorem 77.2]{CR}).

Let $M\neq 0$ be a projective $FG$-module. Then the corresponding
 projective $R G$-module
is called {\it the lifting of} $M$, which we often denote by
$\tilde M$. Obviously, $\dim M$ is equal to the $R $-rank of
$\tilde M$. The latter is equal to the dimension of the
$KG$-module obtained from $ M_{R }$ by the extension of the
coefficient ring to the quotient field $K$ of $R $, so we also
write $\dim M$ for the rank of an $R $-module $M$.

\med
For sake of convenience we record the \f easy observation:

\begin{lemma}\label{ax2}
$(1)$ Let $M=M_1\oplus M_2$ be a
direct sum of $FG$- or $R G$-modules. Then $ C_M(G)=
C_{M_1}(G)\oplus C_{M_2}(G) $.

\med
$(2)$ If $H$ is a subgroup of $G$,   $M$  a projective $FG$-module
with lifting  $L$ 
then $\dim C_M(H)=\dim
C_L(H)$.\end{lemma}

Proof. (1) is trivial. As the restriction of a  projective
$G$-module to $H$ remains projective, it suffices to prove (2)
 for $H=G$. Then (2) is obvious if $M$ is the regular
$FG$-module. By (1), this is implies (2) when $M$ is free, and
hence when  $M$ is projective. \hfill $\Box$

\med
The following is well known (see for instance \cite[p. 52]{Fein}):

\begin{lemma}\label{dp8}
Let $G=H\times N$, the direct product of finite groups $H$ and $N$,  
and let $\Phi,\Psi $ be  PIM's for $H,N$ respectively. Then
$\Phi\otimes \Psi$ is a PIM for $G$, and hence $c_{\Phi\otimes
\Psi}=c_\Phi\cdot c_\Psi$.
\end{lemma}

\med
The \f lemma asserts that, for a projective $G$-module $M$ and a
normal subgroup $N$ of $G$, the  $G/N$-module $C_M(N)$ is
projective. This is a rather general fact, but it does not seem to
be recorded in any standard textbook. The proof below is a
variation of that given in \cite[the proof of Proposition
2.2]{MW}.

\begin{lemma}\label{a2}
Let $H$ be a subgroup of $G$ and   $N$  a normal subgroup of $H$.
Let $M$ be a projective $FG$-module.

\med
$(1)$ $\overline{M}_{H/N}$ is a projective $F(H/N)$-module.

\med
$(2)$ Let  $L$ be the lifting of $ M$. Then $\overline{L}_{H/N}$
is the lifting of $ \overline{M}_{H/N}$.

\med
$(3)$ Let $\chi$ be the character of $L$. Then the truncation
$\overline{\chi}_{H/N}$ is the character of $\overline{L}_{H/N}$.
\end{lemma}

Proof. (1) As $M|_H$ is a projective module, it suffices to prove
the statement for $H=G$.

Assuming $H=G$, suppose first that $M$ is the regular $FH$-module.
 Then $M|_N$ is a free $FN$-module of rank $|H:N|$, and hence
$\dim\overline{FH}=|H:N|$. Let $a:=\sum_{n\in N}n\in FN$. Then the
mapping $h:x\ra xa$ $(x\in FG )$ is an $FH$-module \ho whose
kernel $A$ is spanned by the elements $g(n-1)$ $(n\in N,g\in H)$.
Therefore, $h(FH)=FH/A\cong F(H/N)$. As $xa\in \overline{FH}$ and
$\dim F(H/N)=|H:N|= \dim\overline{FH}$, it follows that
$\overline{FH}$ is isomorphic to  $F(H/N)$, the regular
$F(H/N)$-module.

 Therefore, the lemma is true
if $M$ is free (and in this case $\overline{ M}$ is free).
Otherwise, let $M'$ be a projective $FG$-module such that $M\oplus
M'=J$, where $J$ is free. Then $\overline{ J}$ is a free
$F(G/N)$-module, and $\overline{ M}\oplus \overline{
M'}=\overline{ J}$. As $\overline{ J}$ is a free $F(G/N)$-module,
$\overline{ M}$ is projective.

\med
(2) Let $\pi: L\ra M$ be the reduction map. Then $\pi(\overline{
L})\subseteq \overline{ M}$. By Lemma \ref{ax2}, $\dim\overline{
L}=\dim\overline{ M}$, as desired.

\med
(3) follows from the definition of $\overline{\chi}_{H/N}$. \hfill $\Box$

\med
 For a $p$-group $P$, every projective
$FP$-module $M$ is free, that is, a direct sum of copies of the
regular $FP$-module, which is the only PIM for $P$ (see \cite[\S
54, Exercise 1]{CR} or \cite[Ch.III, Corollaries 2.9, 2.10]{Fe}).
This is therefore true for its lifting as well. \itf  both $M$ and
its lifting have the same rank as $H$-modules, equal exactly to
$\dim M/|P|$. Obviously, $\dim C_M(P)=1$ for the regular
$FP$-module $M$. This implies the following assertion:

\begin{lemma}\label{y2}
Let $P$ be a $p$-subgroup of $G$,
$M$ a projective $FG$-module and $L=M_{K }$. Then
 $ \dim C_M(P)=(\dim M)/|P|=(\dim L)/|P|=\dim C_L(P)$.\end{lemma}

Let  $U\in {\rm Syl}_p(G)$  and $M$  a projective $FG$-module.
Then $c_M:=(\dim M)/|U|$ is an integer. Some formulas become
simpler if one uses $c_M$ instead of $\dim M$. The first equality
of Lemma \ref{y2} implies:

\begin{lemma}\label{x2}
Let $H$ be a subgroup of $G$
with $(|G:H|,p)=1$. Then $c_M=c_{M|_H}$ for any projective
$FG$-module $M$. \end{lemma}

Let $M$ be a projective $FG$-module with  lifting $L$. For
brevity, the character $\chi$ of $L$ is also called the character
of $M$. \itf $\chi$ vanishes at the $p$-singular elements (indeed,
if $g=su$, where $u\neq 1$ is a $p$-element, $s$ is a $p'$-element
and $[s,u]=1$, then $L$ is a direct sum of the eigenspaces of $s$;
by the Krull-Schmidt theorem, every $s$-eigenspace is a free $
R\lan u\ran $-module (as so is $L$), and hence the trace of  $su$
is $0$). See \cite[Ch.IV, Corollary 2.5]{Fe}.

\def\prm{projective module }

\begin{lemma}\label{t1}
Let $G$ be a finite group and
 $U\in {Syl}_p(G)$. Let $M$ be a projective $FG$-module with
lifting $L$, and $\chi$ the character of $L$. Then $c_M=(\chi|_U,
1_U)$. Moreover, $C_M(U)$ and $C_L(U)$ (the fixed point subspaces
of U on $L,M$, resp.)
 are $N_G(U)$-modules with the same Brauer character. \end{lemma}

Proof. As $N_G(U)/U$ is a $p'$-group, $N_G(U)$ splits as $UH$,
where $H\cong N_G(U)/U$. Therefore,  $M|_H$ and $L|_H$ have the
same Brauer characters. Let  $\rho:L\ra M$ be the reduction \hot
Obviously, $\rho(C_L(U))\subset C_M(U)$. As $\dim C_L(U)=c_M=\dim
C_M(U)$, we have $\rho(C_L(U))=C_M(U)$. \itf $C_M(U)|_H$ and
$C_L(U)|_H$ have the same Brauer characters, and the lemma
follows. \hfill $\Box$

\begin{lemma}\label{om1}
Let $N$ be a normal subgroup of
$G$ with $(|G:N|,p)=1$. Let $M$ be an $FG$-module.

\med
$(1)$ $\soc (M|_N)=(\soc M)|_N $.

\med
$(2)$ Let $M$ be a PIM for $G$ and let $S=\soc M$. Suppose that
$S|_N$ is irreducible. Then $M|_N$ is a PIM for $N$.
\end{lemma}

Proof. (1) Let $X$ be an \ir submodule of $M|_N$. Then so is $gX$
for every  $g\in G$, and hence $Y=\sum_{g\in G}gX$ is an
$FG$-module, obviously, completely reducible. Therefore,
$Y\subseteq \soc M$. As $\soc (M|_N)$ is the sum of \ir submodules
of $M|_N$, by the above we have $\soc (M|_N)\subseteq (\soc M)|_N
$. The converse inclusion follows from Clifford's theorem.

(2) By (1), $\soc (M|_H)$ is irreducible, so $M|_H$ is a PIM (as
it is projective). \hfill $\Box$

\med
The following lemma is a special case of \cite[Ch.IV, Lemma
4.26]{Fe}. Our proof below is different, and somehow, more
natural.

\begin{lemma}\label{a4}
 Let $N$ be a normal $p$-subgroup
of $G$ and let $\Phi$ be a PIM for $G$. Then
$\overline{\Phi}:=C_{\Phi}(N)=\soc (\Phi|_N) $ is a PIM for
$G/N$ and $c_\Phi 
=c_{\overline{\Phi}}$. \end{lemma}

\def\hw{highest weight }

 Proof. By Lemma \ref{a2}, $\overline{ \Phi}$ is a
projective $F(G/N)$-module. As $N$ acts trivially on every \ir
$FG$-module, $\soc\Phi\subseteq \overline{ \Phi}$. In fact,
 $\soc\Phi=\soc \overline{ \Phi}$
 since $G$ acts in $\overline{\Phi}$ via $G/N$.
Recall that a projective module
 is a PIM \ii its socle is irreducible. Therefore, $\soc\Phi$ is
 irreducible as an $FG$-module, and hence as an $F(G/N)$-module.
 So $\overline{ \Phi}$ is a PIM for $G/N$, as claimed.

Let $U\in {\rm Syl}_p(G)$ and $\overline{ U}=U/N$.
Obviously, $C_\Phi(U)=C_{\overline{\Phi} }(\overline{ U})$, so
$c_M=\dim C_\Phi(U)=\dim C_{\overline{\Phi}
}(\overline{U})=c_{\overline{\Phi}}$, where $\overline{ \Phi}$ is
viewed as $F(G/U)$-module.  \hfill $\Box$

\begin{lemma}\label{a1}
 Let $H$ be a subgroup of  $G$
and $M$ a projective  $FG$-module. Then $M|_H$ is a projective
module, and $c_{M}= \frac{|H|_p}{|G|_p}\cdot c_{M|_H}$.
\end{lemma}

Proof. The first claim is well known. The second one follows by
dividing by $|G|_p$ the left and the right hand sides of the
equality $|G|_p\cdot c_{M} =\dim M|_H= c_{M|_H}\cdot |H|_p$.
\hfill $\Box$

\begin{lemma}\label{mp1}
Let $H\subset G$ be finite groups such that $H$ contains a \syl of
$G$, and $N=O_p(H)$. Let $M$ be a projective $FG$-module with
socle  $S$. Let  $D=\soc \overline{S}_{H/N}$ and let $M'$ be a
projective $F(H/N)$-module with socle $D$. Then
$\overline{M}_{H/N}$ contains a submodule isomorphic to  $M'$. In
addition, $c_M\geq c_{M'}$.\end{lemma}

Proof. As $S\subset M$ and $N$ is a $p$-group, we have $\soc
(S|_H)\subseteq  \soc (M|_H)\subseteq C_M(N)$. Viewing each module
as an $F(H/N)$-module, we have $D=\soc
(\overline{S}_{H/N})\subseteq \soc (\overline{M}_{H/N})\subseteq
\overline{M}_{H/N} $, and $\overline{M}_{H/N}$ is projective by
Lemma \ref{a2}. Therefore,  $M'\subseteq \overline{M}_{H/N}$. The
additional claim follows from Lemma \ref{a1}, as $|H|_p=|G|_p$.
\hfill $\Box$

\begin{lemma}\label{z1}
Let $H\subset G$ be finite groups. Suppose that $(|G:H|,p)=1$ and
$\Phi $ is a PIM of dimension $|G|_p$, that is,  $c_\Phi=1$. Then
$\Psi:=\Phi|_H$ is a PIM for $H$, and $c_\Psi=1$. In addition,
$\soc (\Phi|_H)$ is irreducible.\end{lemma}

Proof. This follows straightforwardly from Lemma \ref{a1}. \hfill
$\Box$

\begin{lemma}\label{b1}
Let $G$ be a finite group and $U$ a \syl of $G$. Let $\eta:G\ra
\CC$ be a conjugacy class function such that $\eta(g)=0$ for every
$1\neq g\in U$. Let $\tau$ be an \ir character of $U$. Then
$(\eta, \tau^G)=\eta(1)\tau(1)/|U|.$ In particular, if $\eta$ is a
character of a PIM $\Phi$ then $(\eta, \tau^G)=\tau(1)\cdot
c_{\Phi}$.
\end{lemma}

Proof. By Frobenius reciprocity we have: $(\eta, \tau^G)=(\eta|_U,
\tau)=\eta(1)\tau(1)/ |U|$.

\begin{corol} \label{gg1}
Let $G$ be a Chevalley group
in defining characteristic $p$, $U$ a \syl of $G$ and let $\eta$
be the character of a $p$-modular PIM $\Phi$. Then $(\eta,
1_U^G)=c_{\Phi} =(\eta, \Gamma)$, where $\Gamma$ denotes a
Gelfand-Graev character of $G$.
\end{corol}

Note that every Gelfand-Graev character of $G$ is induced from a
certain one-dimensional character of $U$, see \cite{DML} or
\cite{DM}.

\begin{lemma}\label{de4}
Let ${\mathbf H } $ be a reductive algebraic group and $Fr$ a
Frobenius endomorphism of ${\mathbf H } $. Let ${\mathbf G } $ be
the semisimple part of ${\mathbf H } $ and $H={\mathbf H }^{Fr} $,
$G={\mathbf G }^{Fr} $. Then:

$(1)$ Every \ir $FH$-module $M$ remains \ir under restriction to
$G$. Consequently, if $\Phi$ is a PIM for $H$ then  $\Phi|_G$ is a
PIM for $G$.

$(2)$ Let $\Psi$ be a PIM for $G$ with character $\eta$ and let
$\lam\in\Irr G $. For $h\in H$ denote by $\lam^h$ the $h$-twist of
$\lam$. Then $(\lam, \eta)=(\lam^h, \eta)$. In other words, the
rows of the decomposition matrix of  $G$ corresponding to
$H$-twisted ordinary characters coincide. 
\end{lemma}

Proof. (1) It is known that $H$ is a group with BN-pair
\cite[24.10]{MT}.
Let $B$ be a Borel subgroup
of $H$ and $U=O_p(H)$. By \cite[Theorem 4.3(c)]{C1}, $\dim
C_{M}(U)=1$. As $H/G$ is a $p'$-group, $U\subset G$. By Clifford's
theorem, $M|_G$ is completely reducible, and if $M|_G$ is
reducible then $\dim C_{M}(U)>1$, which is false.

The additional statement in (1) follows from Lemma \ref{om1}.

(2) It follows from (1) that $\Psi=\Psi^h$. So  $(\lam,
\eta)=(\lam^h, \eta^h)= (\lam^h, \eta)$.  \hfill $\Box$

\begin{propo}\label{m75}
Let $H\subset G$ be finite groups such that $(|G:H|,p)=1$  and
$N=O_p(H)$. Let $M$ be a projective $FG$-module with socle $S$ and
character $\chi$. Let $D \subseteq\soc (S|_H)$ be an \ir
$F(H/N)$-module and let $\eta$ be the Brauer character of $D$. Let
$\lam\in\Irr H/N$. Suppose that $\eta$ is a constituent of
$\lam\,({\rm mod}\,p)$ with \mult $d$. Then $(\chi, \lam^{\#
G})\geq d$.
\end{propo}

Proof. Let $R$ be the projective $F(H/N)$-module with socle $D$
and character $\rho$. In fact, $R$ is a PIM as $D$ is \irt Then
$(\rho,\eta)=1$ by the Brauer reciprocity \cite[Lemma 3.3]{Fe}.
Furthermore, $(\overline{\chi},\eta)\geq 1$ as $R\subseteq
\overline{M}$. As $\lam\,({\rm mod}\, p)$ contains $\eta$ with
\mult $d$, it follows that $(\overline{\chi},\lam)\geq d$. By the
Harish-Chandra reciprocity \cite[70.1(iii)]{CR2}, we have $(\chi,
\lam^{\# G}) = (\overline{\chi}_{H/N}, \lam)\geq d$. \hfill $\Box$

\begin{corol}\label{xx1}
Let $G$ be a  Chevalley group, and $P$ be a parabolic subgroup.
Let $\Phi$ be a PIM with character $ \chi$ and socle $S$. Suppose
that $\soc (S|_P)$ lifts, and let $\lam $ be  the  character of
this lift. Then $(\chi,\lam ^G)>0$.
\end{corol}

Proof. This is a special case of Proposition \ref{m75} with
$M=\Phi$ and $H=P$. \hfill $\Box$

\section{Lower bounds for PIM dimensions}

Let $G$ be a quasi-simple Chevalley group, and let $B,N $ be
subgroups defining a $BN$-pair structure of $G$. Here $B$ is a
Borel subgroup of $G$, $U=O_p(B)$  and let $T_0$ be a maximal
torus of
$B$.  Then 
 $W_0=N/T_0$ is
 the Weyl group of $G$ as a group with $BN$-pair. If $G$ is
 non-twisted then $W_0$ coincides with $W$, the Weyl group of
 $\GC$.

Every \ir character $\beta$ of $T_0$ inflated to $B$ yields a
character of $B$, trivial on $U$, which we denote by $\beta_B$.
Obviously, $\beta\ra \beta_B$ is a bijection between $\Irr T_0$
and the 1-dimensional characters of $B$ trivial on $U$. Therefore,
the induced \rep $\beta_B^{G}$ coincides with $\beta^{\# G}$.

Recall that a PIM $\Phi$ has an \ir socle $S$;  so the socle of
$\Phi|_B$ contains the socle of $S|_B$.

\begin{lemma} \label{m12}
Let $G$ be a Chevalley group viewed as a group with BN-pair (so
$W_0:=W(T_0)$ is the Weyl group of the BN-pair). Let $\chi$ be a
character of $G$ vanishing at all unipotent elements $g\neq 1$.
Then $c_\chi=(\chi, 1_U^G)= \sum _{\beta }|W_0\beta|\cdot
(\chi,\beta_B^G)=\sum _{\beta }|W_0\beta|\cdot
(\overline{\chi}_{T_0},\beta)$, where $\beta$ runs over
representatives of the $W_0$-orbits in $\Irr T_0$. In particular,
if $\chi$ is the character of a PIM $\Phi$ then $c_\Phi= \sum
_{\beta }|W_0\beta|\cdot (\overline{\chi}_{T_0},\beta).$
\end{lemma}

Proof. By Lemma \ref{b1}, \ref{gg1}, $|U|\cdot (\chi, 1_U^G)=\chi(1)$. Note that
$1_U^G=\oplus _{\beta_0\in \Irr T_0 }\beta_B^G$. Let $\beta'\in \Irr T$. By \cite[Theorem
47]{St}, $\beta_B^G$ and $\beta_B^{\prime G}$ are equivalent \ii
$\beta$ and $ \beta'$ are in the same $W_0$-orbit.
 \itf
$1_U^G=\oplus _{\beta}|W_0:C_{W_0}(\beta)|\beta_B^G$, where
$\beta$ runs over representatives of the $W_0$-orbits in $\Irr
T_0$. This implies the first equality, while the second one
follows from the Harish-Chandra reciprocity formula 
$(\chi,\beta_B^G)=(\overline{\chi}_{T_0},\beta)$. If $\chi$ is the character of
 $\Phi$ then $c_\Phi=\chi(1)/|U|$. \hfill $\Box$

\begin{lemma}\label{2ya}
Let $S$ be a finite $p$-group with normal subgroup $K$. Let $\chi$
be a character vanishing at all non-identity $p$-elements of $S$.
Then $\chi=c_{\chi}\cdot \rho^{reg}_S$ and
$c_\chi=c_{\overline{\chi}_{S/K}}$ (where $\rho^{reg}_S$ denotes
the regular character of $S$).
\end{lemma}

Proof. Let $\tau$ be an \ir character of $S$. Then
$(\chi,\tau)=\tau(1)\chi(1)/|S|=c_\chi\cdot \tau(1)$, whereas
$(\rho^{reg}_S,\tau)=\tau(1)$. So the former claim follows. Let
$M$ be the $\CC S$-module with character $\chi$. \itf $M$ is a
free  $\CC S$-module of rank $c_\chi$. The latter claim is obvious
for the regular $\CC S$-module in place of $M$, which implies the
lemma. \hfill $\Box$

\begin{lemma}\label{x2a}
Let $H$ be a finite group and $U=O_p(H)$. Let $\chi$ be a
character of $G$ vanishing at all non-identity $p$-elements of
$G$. Then $\overline{\chi}_{H/U}$ vanishes at all non-identity
$p$-elements of $H/U$ and
$c_\chi=c_{\overline{\chi}_{H/U}}$.\end{lemma}

Proof. Let $M$ be a $\CC H$-module with character $\chi$, and let
$M'$ be the fixed point subspace of $U$ on M. Note that $x\in M'$
\ii $x=\frac{1}{|U|}\sum_{u\in U} um$ for some $m\in M$. Let $g\in
H$, $u\in U$. Suppose that the projection of  $g$, and hence of
$gu$, into $H/U$ is not a $p'$-element   (the projections of $gu$
and $g$ in $H/U$ coincide).  \itf
$\overline{\chi}_{H/U}(g)=\frac{1}{|U|}\sum_{u\in U}\chi(gu)=0$ by
assumption, whence the first claim. The equality
$c_\chi=c_{\overline{\chi}_{H/U}}$ follows from Lemma
\ref{2ya}.\hfill $\Box$

\begin{propo}
\label{r2b}
Let $G,B=UT_0,W_0$ be as in Proposition $\ref{m12}$, and let
$\chi$ be a character of $G$ vanishing at all $p$-elements $1\neq
u\in G$. Let $\beta\in\Irr T_0$ and $\beta_B$ the inflation of
$\beta$ to $B$. Suppose that $\beta$ is
 an \ir constituent of
$  \overline{\chi}_{B/U}$. Then $c_{ \chi} \geq |W_0\beta|$. In
addition, if $\chi(1)=|G|_p$ then $  \overline{\chi}_{B/U}$ is \ir
and $W_0$-invariant.
\end{propo}

Proof. By \cite[Theorem 47]{St}, $\beta_B^G$ is equivalent to
$w(\beta)_B^G$ for every  $w\in W_0$. By the Frobenius
reciprocity,
$(\chi|_B,\beta_B)=(\chi,\beta_B^G)=(\chi,w(\beta)_B^G)=
(\chi|_B,w(\beta)_B)$. Therefore,   both $\beta_B$ and $
w(\beta)_B$ occur in $\chi|_B$ with equal multiplicity. As
$w(\beta)_B$ is trivial on $U=O_p(B)$, it follows that $w(\beta)$
is a constituent of
  $ \overline{\chi}:=\overline{\chi}_{B/U}$. So  $\overline{\chi}(1)\geq |W_0\beta|$.
  As $B/U$ is a $p'$-group,  $c_{ \overline{\chi}}= \overline{\chi}(1)$. We know that
  $c_\chi\geq c_{ \overline{\chi}}$ (Lemma \ref{2ya}).
  So the result follows. This also implies
 the additional statement, as $1=c_\chi\geq c_{ \overline{\chi}}\geq |W_0\beta|$
 means that $\beta$ is $W_0$-stable.
  \hfill $\Box$

\begin{propo}
\label{r2}
Let $G,B=UT_0,W_0$ be as in Proposition $\ref{m12}$. Let $\Phi$ be
a PIM with socle $S$ and character $\chi$, and  let $\beta_B$ be
the Brauer character of $\soc S|_B$. Then $c_{ \Phi} \geq
|W_0\beta|$.
\end{propo}

Proof. Let $\chi$ be the character of $\Phi$, that is, the
character of the  lifting  $M$ of $\Phi$. By Lemma \ref{t1},
$c_\Phi=(\chi|_U, 1_U)=c_\chi$, and the character of $C_M(U)|_B$
coincides with the Brauer character of $C_\Phi(U)|_B$. Therefore,
$\beta_B$ occurs as a constituent of $C_M(U)|_B$ so
$(\chi|_B,\beta_B)\geq 1$, and hence $\beta$ is a constituent of
$\overline{\chi}_{B/U}$. So the result follows from Proposition
\ref{r2b}. \hfill $\Box$

\med
Remarks. (1) Let $\GC$ be the algebraic group defining $G$ as
$G=\GC^{Fr}$. Then in  Proposition \ref{r2} $S=V_\mu|_G$ for some
\ir $\GC$-module $V_\mu$, where $\mu$ is the highest weight of
$V_\mu$. Moreover, $\beta=\mu|_{T_0}$. If $G$ is a non-twisted
Chevalley group  then $W(T_0)=W(G)$. Therefore, for non-twisted
groups the result coincides with that of Ballard \cite[Corollary
5.4]{Ba}, see also \cite[9.7]{Hub}. (2) Recall that $\beta_B$ is
irreducible, whereas $\overline{\chi}_{B/U}$ may be reducible.

\med
This can be generalized to a parabolic subgroup $P$ in place of a
Borel subgroup $B$, and a Levi subgroup $L$ of $P$ in place of
$T_0$.  However, the statement has to be modified.
 For this we need to replace $W_0$ by a certain group $\overline{W}_L$,
which is contained in $N_G(L)/L$. Specifically, we may assume that
 $ B\subseteq P$, and that $T_0\subseteq L$. (The equality holds
only if $P=B$.) Using the data $B,N,W_0,$ defining the the
$BN$-pair structure of $G$, we set $N_L=\{n\in N: nLn\up =L\}$.
Then $\overline{W}_L=N_L/(N_L\cap L)$. (Note that $\overline{W}_L$
is not the Weyl group of $L$ viewed as a group with $BN$-pair; the
latter is  $(N_L\cap L)/T^0$.) For a character $\lam\in \Irr L$
and $n\in N_L$ one considers the $n$-conjugate $\lam^n$ of $\lam$.
Of course, $\lam^n=\lam$ if $n\in L$. Therefore, $\lam^n$ depends
only on the coset $w:=n\cdot (N\cap L)$, which is an element of
$\overline{W}_L$. So one usually writes $\lam^w$ for $w\in
\overline{W}_L$, with the meaning that $\lam^w=\lam^n$ for $n$
from the pullback of $w$ in $N_L$. If $L=B$ then $\overline{W}_L$
is exactly $W_0$.

\medskip

 Recall   (see Notation)
that $\overline{\chi}_L$  denotes the Harish-Chandra
restriction (or the truncation) of $\chi$, and $\overline{\chi}_L$ coincides with
  $\chi'|_L$, where
$\chi'$ is a character of $P$ trivial on $O_p(P)$ and such that
$\chi|_P=\chi'+\mu$ for some character $\mu$ whose all \ir
constituents are non-trivial on $O_p(P)$.

\begin{lemma}\label{bm1}
Let $P$ be a parabolic subgroup of $G$ and $L$ a Levi subgroup of
$P$. Let $\chi$ be a character of $G$.  Then $\overline{\chi}_L$
is $\overline{W}_L$-invariant. In particular, if $P=B$ then
$\overline{\chi}_{T_0}$ is $W_0$-invariant.
\end{lemma}

Proof. Let $\lam'$ be an \ir \ccc of $ \chi'$, and
$\lam=\lam'|_L$. By the Frobenius reciprocity, $(\chi, \lam^
{\prime ^G})=(\chi|_P, \lam')=(\chi',\lam')$ as $\lam'$ is trivial
on $O_p(P)$. Then $(\chi',\lam')$ equals $(\chi'|_L,
\lam)=(\overline{\chi}_L, \lam)$. Furthermore,
  $\lam^{\prime ^G}=w(\lam')^{\# G}$ for every
$w\in \overline{W}_L$, see \cite[70.11]{CR2}. Hence
$(\chi',\lam)=(\chi',w(\lam))$, and the result follows. \hfill $\Box$

\med
The group $N_L$ acts on $L$ by conjugation, and hence
$\overline{W}_L$ acts on $\Irr L$. Note that for any finite group
$G$ the correspondence $\Phi\ra \soc \Phi$ is compatible with the
\au group action. In other words, if $h$ is an \au of $G$ and
$\Phi^h$ is the $h$-twist of $\Phi$, then $\soc \Phi^h=(\soc
\Phi)^h$.

\begin{propo}
\label{r55}
Let $G$ be a Chevalley group, and let $P$ be a parabolic subgroup
of $G$ with Levi subgroup $L$. Let $\Phi$ be a PIM with socle $S$.
Let $S_1=\soc (S|_P)$   and $S_L=S_1|_L$. Let $\Psi$ be the
projective $FL$-module with socle $S_L$. Then $c_\Phi\geq
|\overline{W}_L:C_{\overline{W}_L}(S_L)|\cdot c_\Psi$.

\end{propo}

Proof. Let $M$ be the lifting of $\Phi. $ Then $\overline{M}_L$
and $\overline{\Phi}_L$ are projective $L$-modules with the same
character $\overline{\chi}_L$, see Lemma \ref{mp1}. (By convention
we call $\overline{\chi}_L$
the character of $\overline{\Phi}_L$.) 
Note that $O_p(L)=1$, so $\overline{\chi}_L$ coincides with the
character
in Lemma \ref{bm1}, which  tells us that $\overline{\chi}_L$ 
 is $\overline{W}_L$-invariant. As a projective $FL$- and $ RL$-module
 is determined by its character, it follows that
$\overline{M}_L^w=\overline{M}_L$ and
$\overline{\Phi}_L^w=\overline{\Phi}_L$ for every $w\in \overline{W}_L$.

By Lemma \ref{sdt}, $S_L$ is \irt As $O_p(P)$ acts trivially on
$S_L$, it follows that $S_L\subseteq \soc \overline{\Phi}_L $.
Then $\Psi\subseteq \overline{\Phi}_L$. By the comment
prior the proposition, the $\overline{W}_L$-orbits of 
$S_L$ and $\Psi$ are of the same size
$l:=|\overline{W}_L:C_{\overline{W}_L}(S_L)|$. As
$\overline{\Phi}_L$ is $\overline{W}_L$-invariant, every $\Psi^w$
($w\in \overline{W}_L$) is in $\overline{\Phi}_L$. Therefore,
$\overline{\Phi}_L$ contains at least  $l$ distinct PIM's
$\overline{\Psi}^w$. Obviously, $\dim \Psi^w=\dim\Psi$ for $w\in
\overline{W}_L$, so $\dim \overline{\Phi}_L\geq l \cdot \dim\Psi$,
and hence $c_{\overline{\Phi}_L}\geq l\cdot c_{\Psi}$. By Lemma
\ref{mp1}, $c_\Phi\geq c_{\overline{\Phi}_L}\geq l\cdot c_{\Psi}$,
as required. \hfill $\Box$

\med
In the remaining part of this section we discuss the question when
the lower bounds provided in Proposition \ref{r2} and Proposition
\ref{r55} are efficient. It is known from Ballard's paper
\cite{Ba} that the bound is sharp for some PIM's, and more
examples are provided in \cite[\S 10.7]{Hub}. However,  in general
the bound is not sharp, and, especially for twisted groups,  there
are some characters $1_{T_0}\neq \beta\in \Irr T_0$ for which the
bound is too small for efficient use. The situation is better for
some non-twisted Chevalley groups; this will be explained in the
rest of this section.

If $C_{W_0}(\beta)=W_0,$ then Proposition \ref{r2} gives
$c_\Phi\geq 1$, which is trivial.
 This always happens if $G$ is a
non-twisted Chevalley group $G(q)$ with $q=2$, or, in general, if
$\beta=1_{T_0}$. In fact, there are more cases where
$C_{W_0}(\beta)=W_0.$ In addition, one needs to decide what is the
minimum size of $|W_0\beta|$ if it is greater than 1. Thus, we are
faced with two problems:

\med
(1) Determine $\beta\in \Irr T_0$ such that $C_{W_0}(\beta)=W_0,$ and

\med
(2) Assuming $|W_0\beta|>1$, find a lower bound for $|W_0\beta|$.


\med
We could obtain a full solution to these problems. However, it
seems that for the purpose of this paper we need only to describe
favourable situations, where the orbit $W_0\beta$ is not too small
for every $\beta\neq 1_{T_0}$. Our results in this line are
exposed in Propositions \ref{dn5}, \ref{cm2} and \ref{r4}, where
the groups considered are non-twisted. To explain our approach, we
therefore assume that $G$ is non-twisted. In this case $W_0$
coincides with the Weyl group $W$ of $\GC$.

Let $\GC$ be a simple simply connected algebraic group in defining
characteristic $p$, and $G=G(q)$. Let $r$ be the rank of $\GC$ and
let ${\bf T_0}$ be a maximal torus of $\GC$. Then the action of
$W=N_\GC({\mathbf T_0})$ on ${\bf T_0}$ yields the action of $W$
on ${\mathbb Z}^r$, the group of rational characters of ${\bf
T_0}$. In turn, this yields a \rep $\zeta_0:W\ra GL(r,{\mathbb
Z})$, which we call  the natural \rep of $W$. It is well known
that $\zeta_0(W)$ is an \ir group (of $GL(n,\CC )$) generated by
reflections. (See \cite[0.31]{DM}.) The \rep $\zeta_0$ is well
understood, see \cite{St1,St}.
 It turns out that, if $\beta\neq 1_{T_0}$, then $|W\beta|$ is
 not too small provided  $\zeta_0(W)$
remains \ir  modulo every prime dividing $q-1$. This requires $q$
to be even for $G$ of type $B,C,D$, see Table 3. For a  prime
$\ell$ dividing $q-1$ denote by $\zeta_\ell$ the \rep obtained
from $\zeta_0$ by  reduction modulo $\ell$. More precise analysis
shows that it is enough that the dual \rep of $\zeta_\ell$, if it
is reducible, were fixed point free. This happens for $SL(r+1,q)$,
$r>1$, when $\ell$ divides $r+1$, and for $E_6(q),E_7(q),E_8(q)$,
see Proposition \ref{cm2}.

\begin{lemma}\label{du3}
Let $G=G(q)$ be a non-twisted Chevalley group of rank $n$, $T_0$ a
split torus, and $W$ the Weyl group of $G$. Let  $\zeta_0:W\ra
GL(n,{\mathbb Z})$ be the natural \rep of $W$ and let $\zeta_\ell$
denote the reduction of $\zeta_0$ modulo a prime $\ell$. Then
$T_0$ has a non-trivial $W$-invariant character \ii there is a
prime $\ell$ dividing $q-1$ such that  $\zeta_\ell$ fixes a
non-zero vector on $F_\ell^n$. In addition, $\zeta_\ell$ is dual
to the natural action of $W$ on $T_\ell$, the subgroup of elements
of order $\ell$ in $T_0$.
\end{lemma}

Proof. If $W $ fixes a character $1_{T_0}\neq \beta\in\Irr T_0$
then it fixes any power $\beta^k$ too. So it suffices to deal with
the case where the order of $\beta$ is a prime. So let
$\ell=|\beta|$ be a prime.

It is well known that $T_0$ is a direct product of cyclic  groups
of order $q-1$, and hence $T_{\ell}=\{t^{(q-1)/\ell}: t\in T_0\}$.
The characters of $T_0$ therefore correspond to elements of
${\mathbb Z}^n/(q-1){\mathbb Z}^n$, and those of  $T_\ell$
correspond to  elements of ${\mathbb Z}^n/\ell{\mathbb Z}^n$. This
yields the reduction mapping $\zeta_0\ra \zeta_\ell$, and the
first assertion of the lemma follows.

The additional claim  describes $\zeta_\ell$ in terms of the
action of $W$ on $T_\ell$. The group $T_0^*:=\Irr T_0$ is
isomorphic to $T_0$, and the actions of $W$ on $T_0$ and on
$T_0^*$ are dual to each other. Let $T_\ell^*=\{t\in
T_0^*:t^\ell=1\}$.
Then
the action of $W$ on $T_\ell^*$ is dual to the action of $W$ on $
T/T^\ell$. As $T$ is homocyclic, the  action of $W$ on $ T/T^\ell$
is equivalent to that on $T_\ell$. \hfill $\Box$


\begin{propo}\label{dn5}
Let $q$ be even, $G=C_n(q)$, $n>1$, or
$D^+_n(q) $,  $n>3$, and let $1_{T_0}\neq \beta\in \Irr T_0$. Then
$|W_0\beta|\geq 2n$.

\end{propo}

Proof. The group $W\cong W(B_n)=W(C_n)$, resp., $ W(D_n)$, is a
semidirect product of a normal 2-group $A$ of order $2^n$, resp.,
$2^{n-1}$, and the symmetric group $S_n$. Note that $|A|\geq 2n$.
In the reflection \rep $\zeta_0$ the group $\zeta_0(W)$ can be
realized by monomial $(n\times n)$-matrices over ${\mathbb Z}$
with diagonal subgroup $A$  and the group $S_n$ as the group of
all basis permuting matrices. This group remains \ir under
reduction modulo any prime $\ell>2$. In addition, $A$ fixes no
non-zero vector on $F_\ell^n$, and if $G=D^+_n(q)$ then $\det a=1$
for $a\in A$.

It suffices to prove the lemma when $|\beta|=\ell$ for every prime
divisor $\ell$ of $|T_0|$. As $|T_0|$ is odd,
  $\ell$  is odd too. So $\zeta_\ell$, the reduction of
  $\zeta_0$ modulo $\ell$,  is an \ir matrix group. Let $0\neq v\in F_\ell^n$. Set
  $X=C_{W}(v)$ and $Y=A\cap X$. Since $A$ acts fixed point freely on ${\mathbb
 Z}^n$,
 and hence on $F_\ell^n$, it follows that $Y\neq A$. We  show that $W:X\geq 2n$.
If $Y=1$ then $|W:X|\geq
  |A| \geq 2n$. Suppose $Y\neq 1$. If $X=Y$ then
$W:X\geq 2|S_n|\geq 2n$.  Suppose  $X\neq
Y$, and let $S=X/Y$. Then $S\subset S_n$.
As   $|A:Y|\geq 2$, we have   $|W:X|= |A:Y|\cdot
|S_n:S|\geq 2\cdot |S_n:S|$. 
If $|A:Y|= 2$ then $v$ is a basis vector, and $S\cong S_{n-1}$. Therefore, $|W:X|=2n$.
Suppose that $|A:Y|>2$. It is easy to check that this implies  that $|Y|=2^i$
with $i>1$ and $X\cong S_i\times S_{n-i}$. This again implies $|W:X|\geq 2n$,
 as required. \hfill $\Box$

\med
In Tables 1, 2, 3  $R$ is an indecomposable root system, and $Z_2$
in Table 2  denotes the cyclic group of order 2. Note that the
data in Tables 2,3 are well known for root systems of types
$A,B,C,D$, and for types $E_i$, $i=6,7,8$, the data follow from
\cite{Atl, MAtl}.

\bigskip
  \begin{center}
 TABLE 1: The structure of  $W(R)$ 
    \vspace{10pt}

{\tiny
   \noindent  \begin{tabular}{|c|c|c|c|c|c|c|c|c| }

        \hline
         $R$&$ A_{n-1} $ &$B_n,C_n$
        &$D_n$&$E_6$&$E_7$&$E_8$&$F_4 $&$ G_2$
             \\
\hline
      $W$&$S_n$ & $ 2^n\cdot S_n$ & $2^{n-1}\cdot S_n $&
      $PSp(4,3)\cdot Z_2 $&$Z_2\cdot Sp(6,2)\cdot Z_2 $&
      $\Omega^+(8,2)\cdot Z_2$&$\Omega^+(4,3)\cdot Z_2$&$Z_2\times S_3$\cr
\hline

    \end{tabular}
}
\end{center}

\bigskip


\bigskip
  \begin{center}
 TABLE 2: The minimum degree of a non-linear \irr of  $W(R)$ 
    \vspace{10pt}

   \noindent  \begin{tabular}{|c|c|c|c|c|c|c|  }

        \hline
         $R$&$ A_{n-1}$,    $B_n,C_n$,
        $D_n$, $n\neq 4$ & $A_3$, $B_4=C_4, D_4$&  $E_6$&$E_7$&$E_8$&$F_4 $, $ G_2$
             \\
\hline
      $d$&$n-1$ & $2$& $ 6$ & $7 $&
      $8 $&$2 $\cr 
\hline

    \end{tabular}
\end{center}


\bigskip
  \begin{center}

 TABLE 3: mod$\,\ell$ irreducibility of the natural \rep of $W(R)$
    \vspace{10pt}

{\small
   \noindent  \begin{tabular}{|c|c|c|c|c|c|c|c|c|c| }

        \hline
        $R$& $A_{n-1} $ &$B_n,C_n$
        &$D_n$&$E_6$&$E_7$&$E_8$&$F_4 $&$ G_2$
             \\
\hline
     $\ell$  & $(\ell,n)\neq 1$ & $ \ell\neq 2$ & $\ell\neq 2 $
     &$\ell\neq 3 $&$\ell\neq 2 $&  any $\ell $&$\ell\neq 2$&  $\ell\neq 3$   \cr
\hline dim &$n-1$&$n$&$n$&$6$&$7$&$8$&$4$&$2$\\ \hline
    \end{tabular}}
\end{center}

\bigskip

\begin{lemma}\label{n1}
Let $d$ be the minimum degree  of a non-linear  \ir character of
$W(R)$. Then $d$ is as in Table $2$.
\end{lemma}

Proof. If $R=A_{n-1}$ then $W(R)\cong S_n$. The degree formula for
\ir \reps of $S_n$ easily implies that    $d\geq n-1$ unless $n=
4$,  where $d=2$. Let $R=B_n,C_n$ or $D_n$. It is well known that
the minimum degree of a faithful \rep of $W(R)$ equals $2n$. As
$W(R)/A\cong S_n $ for an abelian  normal subgroup $A$,  one
arrives at the same conclusion as for $S_n$. If $R=F_4$ then
$W(R)\cong O^+(4,3)$. This group has a normal series $N_1\subset
N_2\subset W(R)$, where $N_1$ is extraspecial 2-subgroup of order
32, $N_2/N_1$ is elementary abelian of order 9, and $W(R)/N_2$ is
elementary abelian of order 4. One observes that $W(R)/N_1$ is
isomorphic to $S_3\times S_3$, and this group has an \ir character
of degree 2. Groups $Sp(4,3)$, $Sp(6,2)$ and $\Omega^+(8,2)$ are
available in \cite{Atl}, so the result follows by inspection.
\hfill $\Box$

\med
Below  $G=E_6(q),E_7(q)$ are groups arising from the simply
connected algebraic group.

\begin{propo}\label{cm2}
Let $G\in \{SL(n,q), n>2, E_6(q),E_7(q),E_8(q), F_4(q), G_2(q)\}$.
Let $W$ be the Weyl group  of  $G$, and $T_0$ a split torus. Then
$1_{T_0}$ is the only $W$-invariant \ir character of $T_0$.

\end{propo}

Proof. We use Lemma \ref{du3} without explicit references to it.

\medskip
Case 1. $G=SL(n,q)$, $n>2$. Here $W\cong S_{n}$. We can assume
that $T_0$ is the group of diagonal matrices. Let $ \zeta_0$ be
the usual \irr of $S_{n}\ra GL(n-1, {\mathbb Z})$. Let $\ell $ be
a prime dividing $|T_0|.$ If $(\ell, n)=1$ then $ \zeta_0$ remains
\ir modulo $\ell$. So assume $(\ell, n)\neq 1$. Then $ \zeta_\ell$
is not completely reducible. It has two composition factors, one
is of dimension $n-2$, and the other factor is trivial (see
\cite[5.3.4]{KL}). Let $\eta$ be a non-trivial $\ell$-root of
unity in $F_q$ and $t=\eta\cdot\Id$; then $t$ is a scalar matrix.
Obviously, $t\in T_0$. It follows that $\zeta^*_\ell(W)$  fixes a
non-zero vector of $F^{n-1}_\ell$. As $\zeta_\ell^*$ is dual to
$\zeta_\ell$, it follows by dimension reason that $\zeta_\ell$ has
no fixed vector, unless, possibly, $n-1= 2$. If $n=3$ then
$\ell=3$. As $\zeta^*_3$ is faithful and reducible, $\zeta^*_3(W)$
is conjugate with the matrix group $\big\{ \begin{pmatrix}1&*\cr 0
&\pm 1\end{pmatrix}\big\}$. Then $\zeta_3(W)$ fixes no non-zero
vector on $F_3^2$.

\medskip
Case 2. $G= E_6(q)$. Here $\zeta_0(W)\subset GL(6, {\mathbb Z})$
and $W\cong PSp(4,3)\cdot Z_2$. Then $\zeta_0$ remains \ir modulo
any prime $\ell\neq 3$, see \cite{MAtl}. Let $\ell=3$, so $3$
divides $q-1$ and $|Z(G)|=(3,q-1)=3$. Therefore,  $W$ fixes a
non-identity element of $T_3$. So $T_3\cong F_3^6$ viewed as a
$W$-module has a one-dimensional subspace $S$, say, fixed by $W$.
As $ PSp(4,3)$ has no non-trivial \irr of degree less than $5$
\cite{MAtl}, it follows that the second composition factor is of
degree 5. As $T_3$ is indecomposable, and the quotient $T_3/S$ is
irreducible, the dual module $T^*_3$ has no trivial submodule.

\medskip
Case 3. $G= E_7(q)$. Here $\zeta_0(W)\subset GL(7, {\mathbb Z})$
and $W\cong Z_2\cdot Sp(6,2)\cdot Z_2$. Then $\zeta_0$ remains \ir
modulo any prime $\ell>2$. Let $\ell=2$. Then $q$ is odd and
$|Z(E_7(q))|=(2,q-1)=2$. Therefore,  the module $T_2$ has a
non-trivial fixed point submodule $S$. As $T_2$ is indecomposable,
and the quotient $T_2/S$ is irreducible, the dual module $T^*_2$
has no trivial submodule.
\medskip

Case 4.
 $G= E_8(q)$. Here $\zeta_0(W)\subset GL(8, {\mathbb Z})$
and $W\cong Z_2\cdot \Omega^+(6,2)\cdot Z_2$. Then  $\zeta_0$
remains \ir modulo any prime $\ell$, see \cite{MAtl}.

\medskip

Case 5.  $G=F_4(q)$. Then $W\cong W(F_4)\cong O^+(4,3)$ and
$\zeta_\ell(W)$ is \ir for $\ell>2$, see Table 2. Note that
$O^+(4,3)'\cong SL(2,3)\circ SL(2,3)$. 
So it suffices to observe that $O^+(4,3)\,$mod$\,2$ fixes no
non-zero vector on $F_2^4$. As this is the case for the Sylow
3-subgroup of $W$,  the result follows.

\medskip

Case 6. $G=G_2(q)$. Then $W$ is the dihedral group of order 12.
Then $\zeta_\ell(W)$ is \ir modulo any prime $\ell\neq 3$.
Let $\ell=3$. Then a Sylow $2$-subgroup of $W$ fixes no
non-non-zero vector in $F_3^2$, and hence this is  true for the
dual action. So again $W$ fixes no element of order 2 of $\Irr
T_0$.  \hfill $\Box$

\begin{propo}
\label{r4}
Let $G\in \{SL(n,q), n>4, E_6(q),E_7(q),E_8(q)\}$. Let
$1_{T_0}\neq \beta\in\Irr T_0$. Then $|W \beta|\geq m$, where
$m=n,  27,28,120$, respectively.
\end{propo}

Proof. By Proposition \ref{cm2}, $|W\beta|>1$.
Let $r$ be the rank of $G$. 
If $W$ is realized via $\zeta_0$ as a subgroup of $GL(r,{\mathbb
Z}) $ generated by reflections then the stabilizer $C_{W}(v)$ of
every vector $v\in {\mathbb Z}^r$ is generated by reflections. Due
to a result of J.-P. Serre, see \cite{KM}, this is also true if
$W$ acts in $F_\ell^r$ and $\ell$ is coprime to $|W|$. This makes
easy the computation of $|W\beta|=|W:C_{W}(v)|$. If $\ell$ divides $|W|$,
and $v\in F_\ell^{r_0}$ then $C_{W}(v)$ is not always generated by
reflections; see \cite{KM} where the authors classify all finite
\ir groups $H$ such that $C_{H}(V)$ is generated by reflections
for every subspace $V$ of $F_\ell^{r}$. Partially we could use
the results of \cite{KM}, but it looks simpler to argue in a more
straightforward way. For our purpose, in most cases it suffices to
know the index of a maximal non-normal subgroup of $W$, which can
be read off from \cite{Atl}.

Let $\ell $ be a prime dividing $|\beta|$. As
$|W:C_{W}(\beta)|\geq |W:C_{W}(\beta^k)|$, it is sufficient to
deal with the case where $|\beta|=\ell$ is a prime.
 Let $D$ be the derived subgroup of
$W$. If
$\zeta_\ell (D)$ 
is irreducible then $|W:C_{W}(\beta)|\geq |D:C_D(\beta)|$, which
is not less than the index $m_D$ of a maximal subgroup  in $D$.

\med
 If $G=SL(n,q)$ 
then $W_0\cong S_{n} $  and  $D\cong {\mathcal A}_{n}$, the
alternating group. It is well known that every proper subgroup of
${\mathcal A}_{n}$, $n>4$, is of index at least $n$. So the lemma
follows in this case. Let  $G=E_6(q)$; then $D\cong SU_4(2)$ and
$m_D=27$, see \cite{Atl}. Let $G=E_7(q)$; then $D\cong Sp_6(2)$
and $m_D=28$ \cite{Atl}.
 If   $G=E_8(q)$ then $D/Z(D)\cong O^+_8(2)'$ and $m_D=120$
\cite{Atl}. \hfill $\Box$

\med
 Remark. If $G=SL(n,q)$,  $q$ odd, $n=4$ then there is
 an element $\beta\neq 1_{T_0}$ with $|W\beta|=3$.
 Indeed, the Sylow 2-subgroup $X$ of $S_4$ has index 3, and
 fixes a non-zero vector of $F_2^3$. By Proposition \ref{cm2},
 this vector is not fixed by $S_4$, whence the claim.
 Let $G=SL(3,q)$. Then $W\cong S_3$, and $\zeta_{\ell}$ is \ir for
 every $\ell\neq 3$. In this case $m=3$.
 Let $\ell=3$. It is observed
 in the proof of Proposition \ref{cm2} that $\zeta_3^*(W)$ fixes no non-zero vector.
 However, the Sylow $3$-subgroup of $W$ fixes a vector $v\neq 0$. Then $|\zeta_3^*(W)v|=2$,
 so $m=2$.

\section{Parabolic descent}

Recall that for a PIM  $\Phi$ of a group $G$
 we set $c_\Phi=\dim\Phi/|G|_p$. If $H$ is a normal subgroup
 of $G$ and $M$ is an $FG$-module, then $C_M(H)$, the fixed point
submodule for $H$, is viewed as $F(G/H)$-module. The socle of a
module  $M$ is denoted by Soc$\,M$. Every PIM is determined by its
socle. The PIM whose socle is $1_G$ is called here 1-PIM, and
denoted by $\Phi_1$.

Let $G$ be a Chevalley group 
so (see Notation)  $G=\GC^{Fr}$, where $\GC$ is simple and simply
connected. Let $P$ be a parabolic subgroup of $G$ with Levi
subgroup $L$. The parabolic descent is the mapping
$\pi_{G,P}:\Phi\ra \Psi$, where $\Phi$ runs over the set of PIM's
for $G$ and $\Psi$ is a PIM for $L$. (One can extend this to the
${\mathbb Z}$-lattices spanned by PIM's for $G$ and $L$.) The
parabolic descent $\pi_{G,P}:\Phi\ra \Psi$ is determined by the
Smith-Dipper mapping $\sigma_{G,L}:\Irr G\ra \Irr L$ defined by
$S\ra \soc \,(S|_{P})$, where $S\in \Irr G$ and the right hand
side is viewed as an $FL$-module.

\begin{lemma}\label{sm9} Let $P$ be a parabolic subgroup of $G$ and $L$ a Levi of $P$.
The mappings $\sigma_{G,P}:\Irr G\ra \Irr L$ and $\pi_{G,P}:{\rm
PIM }_G \ra {\rm PIM }_L$ are surjective.
\end{lemma}

Proof. Let $M$ be an \ir $FP$-module trivial on $O_p(P)$. There
exists an \ir $FG$-module $R$ such that ${\rm Hom}\,(M^G, R)\neq 0
$. By Frobenius reciprocity \cite[10.8]{CR1},  $\dim {\rm
Hom}\,(M^G, R)=\dim {\rm Hom}\,(M, R|_P)$. So $M$ is isomorphic to
a submodule $M',$ say, of  $R|_P$. So $M'\subseteq \soc R|_P$. By
Lemma \ref{sdt}, $\soc R|_P$ is \ir so $M'=\soc R|_P$. This
implies the statement for  $\sigma_{G,P}$. In turn, this  implies
the statement for  $\pi_{G,P}$ as, both for  $G$ and $L$, \ir
modules are in natural bijection with PIM's. \hfill $\Box$

\med
Let $L'$ be the subgroup of $L$ generated by all unipotent
elements of  $L$. If $\Psi$ is a PIM for $L$ then $\Psi|_{L'}$ is
a PIM for $L'$ (Proposition \ref{de4}). Then one may also consider
the mapping $\pi_{G,L'}$ which sends $\Phi$ to $\Psi|_{L'}$. 
For our purpose this version of the parabolic descent has some
advantage. Indeed, there are a parabolic subgroup ${\mathbf P}$ of
${\mathbf G}$ and a Levi subgroup ${\mathbf L}$ of $ {\mathbf P}$
such that   $P={\mathbf P}^{Fr}$ and  $L={\mathbf L}^{Fr}$. Let
${\mathbf L}'$ denote the semisimple part of  ${\mathbf L}$. Then
$L'=({\mathbf L}')^{Fr}$. Thus, $L'$ corresponds to a semisimple
subgroup ${\mathbf L}'$ of $\GC$, and hence the \ir \reps of $L'$
can be parameterized in terms of highest weights. This allows us
to make more precise control of $\pi_{G,L'}$ in terms of
$\sigma_{G,L'}$. By Lemma \ref{mp1}, $c_\Phi\geq c_{\Psi'}$, where
$\Psi'=\pi_{G,L'}(\Phi)$. (This  is useful only if $c_{\Psi'}>1$.)
The main case where $c_{\Psi'}=1$ (and hence $c_\Psi=1$) is when
$\Psi$ is of defect 0. Corollary \ref{fg8} below tells us that if
$\Phi\neq St$ then $\Psi$ is not of defect 0 for some maximal
parabolic subgroup of $G$. In order to prove this, we first turn
Theorem \ref{sdt} to a shape which allows to control $\soc (S|_P)$
in terms of $S$, where $S\in \Irr G$. This is necessary mainly for
twisted Chevalley groups.

\med
  Let $\GC$ be an algebraic group over $F$ of rank $n$,
  $\al_1\ld \al_n$ be simple roots and $\om_1\ld\om_n$ be
the fundamental weights of $\GC$. Let $D$ denote the Dynkin
diagram of $\GC$ with nodes labeled by $1\ld n$ according to
Bourbaki \cite{Bo}. We denote by $X_\al$ the root subgroup of
$\GC$ corresponding to a root $\al$. The dominant weights of $\GC$
are of shape $\sum a_i\om_i$ for some integers $ a_1\ld a_n\geq
0$; for an integer $q$,
 those with $0\leq a_i\leq q-1$ for $i=1\ld n$ are called $q$-restricted.
 The \ir \reps of $\GC$ are parametrized by the dominant weights.
 Given a dominant weight $\mu$, we denote by $V_\mu$ the \irr of $\GC$
 correspondimg to $\mu$; this weight $\mu$ is called the highest weight
 of $V_\mu$.
 For our purpose we may assume that $\GC$ is simply connected.

Let $G=\GC^{Fr}$. Usually one takes for $q$ the common absolute
value of $Fr$ acting on the weight lattice of $\GC$, and set
$G(q)=\GC^{Fr}$. If $q$ is an integer then  $G=G(q)\in
\{SL(n+1,q), SU(n+1,q), Sp(2n,q)$, Spin$\,{}(2n+1,q)$,
Spin$\,{}^\pm(2n,q)$, $E_i(q)$, $i=6,7,8$, ${}^2E_6(q)$, $F_4(q),
G_2(q), {}^3D_4(q) \}$. Otherwise, $q^2$ is an integer, and then
$G\in \{ {}^2B_2(q),{}^2F_4(q), {}^2G_2(q)\}$.

 The \ir \reps of $G$ are parameterized by the dominant weights
 satisfying certain conditions. More precisely,
every \irr of $G$ is the restriction to $G$ of an
\irr of $\GC$ whose highest weight belongs  to the set $\Delta
(G)$, defined   as follows:

$$\Delta (G)=\begin{cases}a_1\om_1+\cdots +a_n\om_n:~0\leq a_1\ld
a_n<q&~ if ~q ~is ~an ~integer,\\ a_i<q\sqrt{1/p}~~ if~~  \al_i~
long, ~and ~a_i<q\sqrt{p}~~ if ~\al_i~ is~ short & ~if ~~q~ is~
not~ an~ integer.\end{cases}$$

\medskip

We refer to the elements of $\Delta (G)$ as  dominant weights for
$G$. (These are called the basic weights for $G$ in
\cite[2.8.1]{GLS}.) By \cite[Theorem 43]{St} $\Delta (G)$
parameterizes the \ir \reps of $G$ up to equivalence. Therefore,
there is a bijection $\Delta (G)\ra \Irr G$, so the \ir \reps of
$G$ can be written as $\phi_\lam$ for $\lam\in\Delta (G)$. Thus,
$\phi_\lam$ extends to a unique \irr of $\GC$ with highest weight
$\lam\in\Delta (G)$, see \cite{St}. For brevity we refer to $\lam$
as the highest weight of $\phi_\lam$.

Furthermore, there a unique weight $\lam\in\Delta (G) $ with
maximal sum $a_1+\cdots +a_n$ (if $q\in\ZZ$ then $a_1=\cdots =
a_n=q-1$). For this $\lam$ $\dim\phi_{\lam}$ is greater than for
all other weights in $\Delta (G)$, and equals $|G|_p$, see
\cite[Corollary of Theorem 46]{St} or  \cite[p. 88]{St1}. The
corresponding $FG$-module is called here {\it the Steinberg
module}, and is denoted by $St$. We record this as follows:

If $q$ is not an integer, then this is refined as follows. Set
$q_1:=q/\sqrt{p}$; then $q_1$  is an integer.

\begin{lemma}\label{st9}
 Define a weight $\si$ as follows: $\sigma=(q-1)(\om_1+\cdots
+\om_n)$ if $q$ is integer, otherwise
and $\sigma=
(q_1-1)\om_1+(2q_1-1)\om_2$, $(3q_1-1)\om_1+(q_1-1)\om_2$,
$(q_1-1)(\om_1+\om_2)+(2q_1-1)(\om_3+\om_4)$, where
$q_1:=q/\sqrt{p}$, respectively, for the group $G={}^2B_2(q)$,
${}^2G_2(q)$, ${}^2F_4(q)$.

Then $\dim V_\si=|G|_p$ and the restriction of $V_{\si}$ to $G$ is
a unique \ir $FG$-module of defect $0$.
\end{lemma}

A standard result of the \rep theory of finite groups implies that
$St$ is a unique \ir $FG$-module of defect 0, and lifts to
characteristic 0. It follows that there is a unique \ir character
of $G$ of degree divisible by $|G|_p$. This is called {\it the
Steinberg character}; usually we keep the notation $St$ for this
character as well.

\med
For a reductive algebraic group $\GC$ Smith's theorem \cite{Sm}
states that if ${\mathbf P}$ is a parabolic subgroup of $\GC$ with
Levi subgroup ${\mathbf L}$ and $V$ is a rational $\GC$-module
then $C_V(O_p({\mathbf P}))$ is an \ir ${\mathbf L}$-module.
Furthermore, suppose that the Frobenius endomorphism stabilizes
${\mathbf P}$ and ${\mathbf L}$, and set $P={\mathbf P}^{Fr}$,
$L={\mathbf L}^{Fr}$. Then $C_V(O_p({\mathbf P}))=C_V(O_p(P))$ and
this is an \ir $FL$-module, see Cabanes \cite[4.2]{Cab}.

To every subset $J\subset D$ one corresponds a parabolic subgroup
${\mathbf P}_J$ by the condition $X_{\al_i}\in {\mathbf P}_J$ for
$i\in D$ and $X_{-\al_i}\in {\mathbf P}_J$ \ii $i\in J$. (These $
{\mathbf P}_J$ are called standard parabolic subgroups. If $J$ is
empty, $ {\mathbf P}_J$ is a Borel subgroup.) Note that for a
subset $J'\in D$ the inclusion $ {\mathbf P}_J\subset {\mathbf
P}_{J'}$ holds \ii $J\subset J'$; in particular, every ${\mathbf
P}_J$ contains the standard Borel subgroup. Set $\GC_J=\lan X_{\pm
\al_i}:i\in J\ran$. Then $\GC_J$ is the semisimple component of a
Levi subgroup ${\mathbf L}_J$ of ${\mathbf P}_J$. If ${\mathbf
P}_J$ and ${\mathbf L}_J$ are $Fr$-stable, then so is ${\mathbf
G}_J$. We set $P_J={\mathbf P}_J^{Fr}$, $L_J={\mathbf L}_J$ and
$G_J=\GC_J^{Fr}$; these $P_J$ are called standard parabolic
subgroups of $G$. 

The \f is known but we have no explicit reference:

\begin{lemma}\label{le1}
$C_V(O_p(P_J))$ is an \ir
$FG_J$-module.\end{lemma}

Proof. By Lemma \ref{sdt}, $C_V(O_p(P_J))$ is an \ir
$FL_J$-module, so the claim follows from Lemma \ref{de4}. \hfill
$\Box$

\med
 There is some
advantage of dealing with ${\mathbf G}_J$ in place of ${\mathbf
L}_J$. The \f result is well known \cite{Sm}:

\begin{lemma}\label{sd1}
  Let $\GC$ be a simple algebraic group over $F$, 
$J$ a non-empty set of nodes at the Dynkin diagram of $\GC$, and
$\GC_J = \langle X_{\pm \al_j} : j\in J\rangle $. Let $V$ be an
\ir $\GC$-module of highest weight $ \om$, and let $v \in V$ be a
vector of  weight $\om$. Then $V_J := \lan \GC_J \, v\ran_{F}$ is
an irreducible direct summand of $V|_{\GC_J}$, with highest weight
$\om_{J}=\sum_{j\in J} a_j\om_j$.
\end{lemma}

If $J$ is connected then $\GC_J$ is a simple algebraic group of
rank  $|J|$, and one may think of the fundamental weights of
$\GC_J$ as $\{\om_j:j\in J\}$. Then $\om_{J}$ means $\sum_{j\in J}
a_j\om_j$. If $J$ is not connected, let $J=J_1\cup\cdots\cup J_k$,
where $J_1\ld J_k$ are the connected components of $J$. Then
$\GC_J$ is the central  product of simple algebraic groups
$\GC_{J_1}\ld \GC_{J_k} $, where $\GC_{J_i}$ corresponds to $J_i$
($i=1\ld k$). Then it is convenient to us to view $\om_J$ as the
string $(\om_{J_1}\ld \om_{J_k})$. Furthermore, $V_J|_{\GC_J}$ is
the tensor product of the \ir \reps of $\GC_{J_i}$ with highest
weight $\om_{J_i}$ for $i=1\ld k$.

\med

There is a version of Lemma \ref{sd1} for finite Chevalley groups.
Lemma \ref{sdt} is insufficient as it does not tell us   how
$(\overline{V}_L)|_{G_J}$ depends on $V$ (in notation of  Lemma
\ref{sd1}). If $G=G(q)$ is non-twisted then this is easy to
describe. Indeed, every \ir $FG$-module extends to a $\GC$-module
with $q$-restricted highest weight; call it $V$.  Then
$G_J:=\GC_J^{Fr}$ is a non-twisted Chevalley group corresponding
to $\GC_J$, and the weight $\om_{J}$ is $q$-restricted. Therefore,
an \ir $\GC_J$-module $V_J$ remains \ir as an $FG_J$-module, and
can be labeled by $\om_J$. In addition, $G_J$ is the central
product of $G_{J_i}:=\GC_{J_i}^{Fr}$.

This argument can be adjusted to obtain a version for twisted
Chevalley group but the twisted group case is less
straightforward. The matter is that $Fr$ induces a permutation
$f$, say, of the nodes of the Dynkin diagram of $G$, which is
trivial \ii $G$ is non-twisted. In the twisted case a set $J$ is
required to be $f$-stable. If every connected component of $J$ is
$f$-stable, then $\om\in \Delta(G)$ implies $\om_J\in
\Delta(G_J)$. So again we can use $\om_J$ to identify
$(V_J)|_{G_J}$.

An additional refinement is required if there is a connected
component $J_1$, say, of $J$ such that $J_2:=f(J_1)\neq J_1$. If
there are roots of different length then, by reordering $J_1,J_2$
we assume that the roots $\al_i$ with $i\in J_1$ are long. Note
that the non-trivial $f$-orbits on $\{1\ld n\}$ are of size $a=2$,
except for the case $G={}^3D_4(q)$ where $a=3$ \cite{St}. Then
$G_J:=(\GC_{J_1}\circ \GC_{J_2})^{Fr}\cong G_{J_1}(q^2)\cong
\GC_{J_1}^{Fr^2} $ if $a=2$, or $G_J:=(\GC_{J_1}\circ
\GC_{J_2}\circ \GC_{J_3})^{Fr}\cong G_{J_1}(q^3)\cong
\GC_{J_1}^{Fr^3}$ if $a=3$. Thus, in this case the Chevalley group
obtained from the $f$-orbit on $J$ is non-twisted and
quasi-simple. So one would wish to identify the \rep $V_J|_{G_J}$
in terms of  algebraic group weights of $\GC_{J_1}$ rather than of
$\GC_{J_1}\circ \GC_{J_2}$ when $a=2$, or $\GC_{J_1}\circ
\GC_{J_2}\circ \GC_{J_3}$ when $a=3$.  We do this in the following
proposition. For this purpose it suffices to assume that $f$ is
transitive on the connected components of $J$. To simplify the
language, we call the highest weight of $\GC_{J_1}$ obtained in
this way the highest weight of $V_J$.

\begin{propo}\label{sd3}
Let $V$ be an \ir $\GC$-module of highest weight $\om=\sum
a_i\om_i$ such that $V|_{G}$ is \ir (so $\om\in\Delta(G)$).  Let
$J$  be an $f$-stable set of nodes at $D$, the Dynkin diagram of
$\GC$. Suppose that $J$ is not connected and $f$ is transitive on
the connected components of $J$. Set $J_i=f^{i-1}(J_1)$ for
$1<i\leq a$ and $\om_{J_i}=\sum_{j\in J_i}a_j\om_j$.  Then
$G_J\cong G_{J_1}(q^a)$.

Let $\tilde{\om}_J$ be the highest weight of $V_J$ viewed as a
$G_{J_1}(q^a)$-module. Then $\tilde{\om}_J\in
\Delta(G_{J_1}^{Fr^(q^a)})$. More precisely, set
$\om'_{J_1}=\sum_{j\in J_1}a_{f(j)}\om _{j}$ and, if $a=3$ set
$\om''_{J_1}=\sum_{j\in J_1}a_{f^2(j)}\om _{j}$. Then:

If $q$ is an integer then $\tilde{\om}_J=\om_{J_1}+q \om'_{J_1}$
for $a=2$, and
$\tilde{\om}_J=\om_{J_1}+q\om'_{f(J_1)}+q^2\om''_{f^{2} (J_1)}$
for $a=3$;

If $q$ is not an integer then  $q^2=p^{2e+1}$ for some integer
$e\geq 0$, and 
$\tilde{\om}_J=\om_{J_1}+ p^{e}\om_{f (J_1)}$.
\end{propo}

Proof. We consider only $a=2$, as the case $a=3$ differs
 only on notation. Thus, we show that $G_J:=(\GC_{J_1}\circ
 \GC_{J_2})^{Fr}\cong G_{J_1}(q^2)$. Note that
 $Fr$ permutes $\GC_{J_1}$ and $ \GC_{J_2}$, and
 acts as follows. Let $x_i\in \GC_{J_i}$ $(i=1,2)$.
 If $q=p^e$ is an integer then $Fr(x_1,x_2)=(Fr^e_0 x_2, Fr^e_0x_1)$,
 where $Fr_0$ is the standard
Frobenius endomorphism arising from the mapping $y\ra y^p$
 $(y\in F)$. If $q$
 is not an integer, then
 $Fr(x_1,x_2)=(Fr^{e+1}_0 x_2, Fr^{e}_0x_1)$.

 So $Fr^{2}$ stabilizes each $\GC_{J_i}$, and its fixed point
 subgroup on $\GC_{J_i}$ is $G_{J_i}(q^2)$. Then $(x_1,x_2)$ is
 fixed by $Fr$ \ii $x_1\in G_{J_1}(q^2)$ and $x_2=Fr_0^e(x_1)$.
 So $(\GC_{J_1}\circ
 \GC_{J_2})^{Fr}\cong G_{J_1}(q^2)$, as claimed. Furthermore,
$(V_J)|_{\GC_{J_1}\circ
 \GC_{J_2}}$ is the tensor product of the \ir $\GC_{J_1}$- and
 $\GC_{J_2}$-modules of highest weights $\om_{J_1}$ and $\om_{J_2}$,
 respectively (as the groups $\GC_{J_1}$ and
$ \GC_{J_2}$  commute elementwise). One can consider the
$\GC_{J_1}$-module obtained from $W|_{\GC_{J_1}\circ
 \GC_{J_2}}$ via the \ho $\GC_{J_1}\ra \GC_{J_1}\circ
 \GC_{J_2}$ defined by $x_1\ra (x_1,Fr_0^e(x_1))$.
  Clearly, this is
 the tensor product of the $\GC_{J_1}$-modules of highest weights
 $\om_{J_1}$ and $p^e\om'_{f(J_1)}$.

If $q\in {\mathbb Z}$ then
 $\om_{J_1}$ and $\om'_{ f(J_1)}$
   are $q$-restricted, so
 $\tilde{\om}_J=\om_{J_1}+q\om'_{f  (J_1)}$ is $q^2$-restricted, and hence
 belongs to  $\Delta(G_{J_1}(q^2))$.

 Suppose that $q\notin {\mathbb Z}$. As $V|_G$ is irreducible,
 $\om\in\Delta(G)$. This implies that $a_i<p^e$ for $i\in J_1$
 and $a_{f(i)}<p^{e+1}$. Then $a_i+p^{e}a_{f(i)}<q^2$, as
 required. (So the highest weight of the $\GC_{J_1}$-module
 in question is
 $\om_{J_1}+p^{e}\om'_{f ( J_1)}\in \Delta(G_{J_1}(q^2))$.) \hfill $\Box$

\med
Remark. In Proposition \ref{sd3} $J$ is disconnected, which
implies that the BN-pair rank of $G$ is at least 2.

\med
Example. Let $\GC$ be of type $A_{2n-1+k}$, $k=0,1$,
$G=SU_{2n+k}(q)$ and $J=\{1\ld n-1,n+k+1\ld 2n-1+k\}$. Then
$J_1=\{1\ld n-1\}$ and $G_{J}\cong SL(n,q^2)$. Let
$\mu=a_1\om_1+\cdots +a_{2n-1+k}\om_{2n-1+k}$. Then
$\tilde{\om}_J=(a_1+qa_{2n-1+k})\om_1'+\cdots
+(a_{n-1}+qa_{n+1+k})\om'_{n-1}$, where $\om'_{1}\ld \om'_{n-1}$
are the fundamental weights of $\GC_{J_1}=A_{n-1}$.

\med
As above,  $D$ denotes  the Dynkin diagram of $\GC$ and
$G=\GC^{Fr}$. Recall that the standard parabolic subgroups $P_J$
of $G$ are in bijection with $f$-stable subsets $J$ of the nodes
of $D$ (and $P_J=G$ \ii $J=D$).

\begin{lemma}\label{st0}
Let $D=J\cup J'$ be the disjoint  union of two $f$-stable subsets.
Then $\Irr G\ra \Irr (G_J\times  G_{J'} )$ is a bijection. In
addition, if $T_J,T_{J'}$ are maximal split tori in $G_J,G_{J'}$,
respectively, then $T_JT_{J'}$ is a maximal split torus in $G$.
\end{lemma}

Proof. The first statement follows from the reasoning in \cite[p.
79]{GLS}. The second one follows from \cite[Theorem
2.4.7(a)]{GLS}, which tells us that $\Pi_i T_{I(i)}$ is a maximal
split torus in $G$, where $I(i)$ is the $f$-orbit containing $i$
and $T_{I(i)}$ is a a maximal split torus in $ G_{I(i)})$.

\med
By induction, it follows that a similar statement is true for any
disjoint union $D=\cup_iJ_i$ of $f$-stable subsets $J_i$ of $D$,
in particular, when  $J_i$ is an $f$-orbit for every $i$. (Only
this case is explicitly mentioned in \cite{GLS}.) Note that the
bijection in Lemma \ref{st0} yields a bijection $PIM_G\ra
PIM_{G_J\times  G_{J'}}$.

\begin{corol}\label{fg8}
$G=\GC^{Fr}$ be a Chevalley group of $BN$-pair rank at least $2$,
and $G\neq {}^2F_4(q)$.  Let $V$ be an \ir $\GC$-module with
highest weight $\mu=\sum a_i\om_i\in \Delta (G)$. Suppose that
$\emptyset\neq J\subseteq D$ is $f$-stable. Let $P:=P_J$ be the
standard parabolic subgroup of $G$ corresponding to $J$, and
$L:=L_J$  a Levi subgroup of $P$. 
Let $V_J=C_V(O_p(P))$.

\med
$(1)$ $\dim V_J=1$ \ii  $a_i=0$ for all $i\in J$.

\med
$(2)$ $V_J|_L$ is of
defect $0$ 
\ii  $a_i=q-1$ for all $i\in J$.
(The latter is equivalent to saying that 
$M|_{G_J}$ is Steinberg).

\med $(3)$ Let $J'=D\setminus J$ and $L'$ be a Levi subgroup of
$P_{J'}$. If $V_J|_{L}$ and $V_{J'}|_{L'}$ are of defect $0$  then
so is $V|_G$. Equivalently, If $V_J|_{G_J}$ and $V_J|_{G_{J'}}$
are Steinberg modules   then $V|_G$ is Steinberg.

\end{corol} \def\ps{parabolic subgroup }
\def\hw{highest weight }

Proof.  Let ${\mathbf P}_J$ 
and 
$\GC_J$ be the respective algebraic groups. Then
$V_J|_{G_J}=C_V(O_p({\mathbf P}_J))$, and $V_J$ is an \ir
$\GC_J$-module.

(1) follows from Lemma \ref{sd1}.

(2) Suppose that $a_i=q-1$ for all $i\in J$. Then $V_J|_{G_J}$  is
Steinberg by Lemma  \ref{st9}. To prove the converse, observe that
if $a_i<q-1$ for some $i\in J$ then $\dim V_J<\dim St$ by Lemma
\ref{st9} applied to ${\mathbf G}_J$. As  $|L|_p=|G_J|_p$, it
follows that $V_J$ is of defect 0. This implies (2).

(3) By (2), $a_i=q-1$ for all $i\in D$, so the claim follows from
Lemma \ref{st9}. \hfill $\Box$

\begin{lemma}\label{st6}
Let $M$ be the Steinberg $FG$-module
for a Chevalley group  $G=\GC^{Fr}$, and let $B$ be a Borel
subgroup of  $G$.  Then $\soc (M|_B)=1_B$. \end{lemma}

Proof. Let $U$ be the \syl  of $B$. As $M$ is projective and $\dim
M=|U|$, it follows that $M|_U$ is the regular  $FU$-module.
Therefore, $\dim\soc (M|_B)=1$. Let $\lam$ be the Brauer character
of $\soc (M|_B)$. As $U$ is in the kernel of $\lam$, it can be
viewed as an ordinary character of  $B$. Let $St$ denote the
character of $M$, so this is exactly the Steinberg character. By
Proposition \ref{m75}, $(St,\lam^G)>0$. If $\lam\neq 1_B$ then
$(\lam^G,1_B^G)=0$ \cite[Theorem 70.15A]{CR2}. In addition,
$(St,1_B^G)=1$ \cite[Theorem 67.10]{CR2}. Therefore, $\lam=1_B$,
as required. \hfill $\Box$

\begin{lemma}\label{sx6}
Let $q\in\ZZ$ and let $V_\mu$ be an \ir $\GC$-module of highest
weight $\mu=\sum a_i\om_i\in\Delta(G)$. 
Let $B$ be a Borel subgroup of  $G=\GC^{Fr}$. Then the \f are
equivalent:

\med
$(1)$ $\soc (V_\mu|_B)=1_B$.

\med
$(2)$ $a_1\ld a_n\in \{0,q-1\}\ $ and $a_{f(i)}=a_i$ for every
$i$.
  \end{lemma}

Proof. Note that the result is well known for $G=SL(2,q)$.

Suppose first that $G\cong SU(3,q)$, so $ \GC\cong SL(3,F)$. Let
$M$ be the \ir module for $\GC$ of highest weight $\om_1$. Then
$M|_{G}$ is isomorphic to the natural module $M'$, say, for
$G\cong SU(3,q)$. Let ${\mathbf B} $ be a Borel subgroup of $\GC$
such that $B={\mathbf B}\cap G$ is a Borel subgroup of $G$. As
explained in the discussion after Lemma \ref{st9}, the socle of
${\mathbf B}$ on $M$ coincides with the socle of $B$. Let $v\in M$
be a vector of highest weight  $\om_1$ on $M$. Then $\lan v\ran$
is ${\mathbf B}$-stable and hence $B$-stable. We can view $M'$ as
an $F_{q^2}G$-module endowed by a unitary form, and choose  $v\in
M'$. Then $v$ is an isotropic vector. It is easy to see that $T_0$
is isomorphic to the multiplicative group $F^{\times}_{q^2}$ of
$F_{q^2}$. Let $\nu\in \Irr T_0$ be the \rep of $B$ on $\lan
v\ran$. Then it is faithful. As $M^*$, the dual of $M$, has
highest weight $\om_2$,  one can check that the \rep of $T_0$ in
$C_{M^*}(B)$ is $\nu^q$. Therefore, the \rep of $T_0$ on
$C_{V_\mu}(B)$ is $\nu^{a_1}\nu^{a_2q}$. \itf $\soc
(V_\mu|_B)=1_B$ \ii $\nu^{a_1}\nu^{a_2q}=1_{T_0}$.   Let $t$ be a
generator of $T_0$. Then $\nu^{a_1}\nu^{a_2q}(t)=t^{a_1+a_2q}$. It
is clear that $t^{a_1+a_2q}=1$ \ii either $a_1=a_2=0$ or
$a_1=a_2=q-1$. This completes the proof for $G=SU(3,q)$.

Other groups $G$ are of BN-pair rank 2, and hence satisfy the
assumptions of Corollary \ref{fg8}.

$(1)\ra (2)$. Let $i\in D$ be such that $a_i\neq 0,q-1$, and let
$J$ be the $f$-orbit of $i$. Let $P_J$ and $G_J$ be as in
Corollary \ref{fg8}, and let $ B_J=G_J\cap B$, so $B_J$ is a Borel
subgroup of $G_J$.

Suppose first that $|J|=1$. Then   $J=\{i\}$ for some  $i$, and
$G_J=SL_2(q)$ (both in the twisted and non-twisted cases). Then
$V_J:=(\soc (V_\mu |_{P_J})|_{G_J}$ is \ir with highest weight
$a_i\om_i$.  Then $\soc (V_J|_{B_J})=1_{B_J}$ \ii
$a_i\in\{0,q-1\}$. As $1_B=(\soc (V_\mu|_B))$ equals $\soc (V_J
|_{B_J})$ inflated to $B$, it follows that $a_i\in\{0,q-1\}$.

Let $f(i)=j$, and $i,j$ are not adjacent. If $|J|=2$ then
$G_J=SL_2(q^2)$ and the highest weight of $V_J$ as an
$FG_J$-module is $a_i+qa_j$. So $a_i+qa_j=q^2-1$ or $0$.
Therefore, $a_i=a_j\in\{0,q-1\}$. Similarly, if $|J|=3$ then
$G_J=SL_2(q^3)$, so $a_i+qa_{f(i)}+q^2a_{f^2(i)}=q^3-1$ or $0$.
This again implies $a_i=a_{f(i)}=a_{f^2(i)}\in\{0,q-1 \}$.

Let $f(i)=j$, $f(j)=i$ and $i,j$ are  adjacent. Then  $G_J\cong
SU_3(q) $ as  $q\in \ZZ$. As $\soc (V_J|_{B_J})\subseteq  \soc
(V_\mu|_B)=1_B$, it follows by the above that  $a_i=a_{f(i)}\in
\{0,q-1\}$.

\med

$(2)\ra (1)$. Let $J=\{i: a_i=q-1\}$ and $J'=\{i: a_i=0\}$. Let
$T_0$ be a maximal torus of $B$.
 Then $\soc (V_\mu|_B)=\soc (\soc (V_\mu|_P))|_B)$
as $\dim \soc (V_\mu|_B)=1$ and $\soc (V_\mu|_P)$ is \irt Let
$B_J$ be a Borel subgroup of $G_J$ and $T_J=T_0\cap B_J$. By
Corollary \ref{fg8}(2), $(\soc (V_\mu|_P)|_{G_J}$ is the Steinberg
$FG_J$-module, so, by Lemma \ref{st6}, $((\soc
(V_\mu|_P)|_{G_J})|_{B_J}=1_{B_j}$, and hence $((\soc
(V_\mu|_P)|_{B_J})|_{T_J}=1_{T_j}$.

Set  $T_{J'}=T_0\cap G_{J'}$. Then $T_0=T_JT_{J'}$, see Lemma
\ref{st0}. 

Let $B_{J'}$ be a Borel subgroup of $G_{J'}\subset P_{J'}$
containing $T_{J'}$. Then $S:=(\soc (V_\mu|_{P_{J'}})|_{G_{J '}}$
is the trivial $FG_{J'}$-module, so
 $S|_{B_{J'}}=1_{B_{J'}}$, and hence $S|_{T_{J'}}=1_{T_{J'}}$.
Therefore, $(\soc (V_\mu|_{B})|_{T_{J'}}=1_{T_{J'}}$. \itf $(\soc
(V_\mu|_{B})|_{T_0}=1_{T_0}$, and then $(\soc (V_\mu|_{B})=1_B$,
as required. \hfill $\Box$

\begin{propo}\label{yg8}
Let  $V_\mu$ be an \ir $\GC$-module of highest weight $\mu=\sum
a_i\om_i\in\Delta(G)$, and let $J=\{i: a_i=q-1\}$. Suppose that
$q$ is an integer, $J$ is $f$-stable and $a_i\in\{0,q-1\}$ for
$i=1\ld n$. Let $B$ be a Borel subgroup of $G$, and let $P=P_J$ be
a parabolic subgroup of $G$. Let $\Phi$ be the PIM with socle
$V:=V_\mu|_G$ and character $\chi$.

\med
$(1)$ $\soc (V_\mu|_B)=1_B$.

\med
$(2)$  $(\chi,1_B^G)>0$.

\med
$(3)$ Let $L$ be a Levi subgroup of $P$ and $S:=\pi_{G,P}(V)$ (which is  an $FL$-module). Then $S$ is of defect $0$. 
Furthermore, let $\rho$ be the character of the lift of $S$.
Then $\rho$ is a constituent of $1_{B\cap L}^{L}$ and $(\chi, \rho
^{\#G})>0$.
\end{propo}

Proof. (1) is contained in Lemma \ref{sx6}. (2) This follows from
Corollary \ref{xx1}.

 (3)  By Corollary
\ref{fg8}, the $S|_{G_J}$ is the Steinberg module. So $S$ is an
$FL$-module of defect 0, and hence lifts to characteristic 0. So
$(\chi, \rho ^{\#G})>0$ by Corollary \ref{xx1}.

Furthermore, by the reasoning in Lemma \ref{sx6},
 $\soc (\soc (V_\mu|_P)|_B)$  coincides with $\soc (V_\mu|_B)=1_B$.
 As $S=(\soc (V_\mu|_P))|_L$, it follows that $S|_{L\cap B}=1_{L\cap B}$.
By Lemma \ref{a2}, where one takes $M=S$, $H=B\cap L$ and
$N=O_p(H)$, the truncation $\overline{\rho}$ is the trivial
character. By the Harish-Chandra reciprocity, $\rho$ is a
constituent of $1_{B\cap L}^{L}$. \hfill $\Box$

\med
The \f lemma will be used in Section 7 in order to determine the PIM's for $G$
of dimension $|G|_p$.

\begin{lemma}\label{sd7}
Let $G=\GC^{Fr}$ be a Chevalley group
of $BN$-pair rank at least $2$, and if $G= {}^2F_4(q)$ assume
$q^2>2$.  Let $n$ be the rank of $\GC$, and let $V$ be an \ir
$F\GC$-module of highest weight $0\neq \mu=a_1\om_1+\cdots
+a_n\om_n\in \Delta(G) $. Suppose that $V|_G\neq St$ and  for
every parabolic subgroup $P$ of $G$ 
the restriction of $\soc (V|_P)$ to a Levi subgroup $L$ is
either projective, or 
 the semisimple part $L'$ of
$L$ is of type $A_1(p)$, $A_2(2)$ or ${}^2A_2(2)$ and $\soc
(V|_P)$ is trivial on $L'$.
Then one of the \f holds:

\med
 $(1)$ $G=SL(3,p)$ and $\om=(p-1)\om_1$ or $(p-1)\om_2;$

\med
 $(2)$ $G=Sp(4,p)\cong {\rm Spin}\,(5,p)$
 and $\om=(p-1)\om_1$ or $(p-1)\om_2;$

\med
 $(3)$ $G=G_2(p)$ and $\om=(p-1)\om_1$ or $(p-1)\om_2;$

\med
 $(4)$ $G=SU(4,p)$ and $\om=(p-1)(\om_1+\om_3);$

\med
 $(5)$ $G={}^3D_4(p)$ and $\om=(p-1)(\om_1+\om_3+\om_4);$

\med
 $(6)$ $G=SU(5,2)$ and $\om=\om_1+\om_4;$


\noindent In addition, in all these cases $\soc (V|_B)=1_B$.
\end{lemma}

Proof. As above, for an $f$-stable subset of nodes of the Dynkin
diagram of $\GC$ we denote by $P_J$ the corresponding parabolic
subgroup of $G$, by $G_J$ the standard semisimple subgroup of
$P_J$ and set $V_J=C_V(O_p(P_J))$ viewed as an $FG_J$-module.

Suppose first that $G$ is of type ${}^2F_4(q)$, $q^2>2$. Then
there are two $f$-orbits on $D$: $J=\{1,4\}$ and $J'=\{2,3\}$.
Accordingly, $G$ has two parabolic subgroups $P_J$, $P_{J'}$. By
Proposition \ref{sd3}, the highest weight of $V_J$ is
$(a_1+2^ea_4)\om'_1$, where $\om'_1$ is the fundamental weight for
$\GC_J$ of type $A_1$ and $2^e=q/\sqrt{2}$. The highest weight of
$V_{J'}$ is $a_2\om'_1+a_3\om'_2$, where $\om'_1$, $\om'_2$ are
the fundamental weights for $\GC_{J'}$ of type $B_2$.

If $V_{(a_1+2^ea_4)\om'_1}$ is the Steinberg module for $G_{J'}$
then $a_1+2a_4=q^2-1$, whence we have
 $a_1=q/\sqrt{2}-1$, $a_4=q\sqrt{2}-1$ (as $0\leq
a_1<q/\sqrt{2}$, $0\leq a_4<q\sqrt{2}$). If 
the restriction of $V_{a_2\om'_1+a_3\om'_2}$ to $ G_{J'}$ is of
defect 0 then $V_{a_2\om'_1+a_3\om'_2}|_{ G_{J'}}$ is  the
Steinberg module, and hence  $a_2=q/\sqrt{2}-1$ and
$a_3=(2q/\sqrt{2})-1$ by Lemma \ref{st9}. Thus,
$a_1=a_2=q/\sqrt{2}-1$, and $a_3=a_4=(2q/\sqrt{2})-1$. By Lemma
\ref{st9}, $V |_G$ is the Steinberg module.

\med
So assume that $G\neq {}^2F_4(q)$. By Lemma \ref{st9}, at least
one of $a_1\ld a_n$ differs from $q-1$.

Suppose first that $G=G(q)$ is non-twisted of rank 2. In notation
of Lemma \ref{sd1} and comments following it, $G_J\cong SL(2,q)$
for $J=\{1\},\{2\}$, so $q=p$, and all these groups are listed
  in items $(1),(2), (3)$. In addition, $V_J$ is
of \hw $a_1\om_1$ or $a_2\om_1$. So the claim about the weights in
  $(1),(2), (3)$ follows.

Suppose that $G=G(q)$ is non-twisted of rank  greater than 2, and
let $D$ be the Dynkin  diagram of $\GC$. Then
 one can remove a suitable edge node  from $D$
 such that $a_i\neq q-1$ for some $i$ in the remaining  set  $J$ of
 nodes.
Then the $FG$-module $V_J$ is not projective (Corollary
\ref{fg8}), and hence we have a contradiction, unless $G_J\in
\{A_1(p),A_2(2)\}$ and $a_i=0$ for $i\in J$. In fact, $G_J\neq
A_1(p)$ as otherwise $|D|=2$. If $G_J\cong A_2(2)\cong SL(3,2)$
then $|D|=3$ and $G=SL(4,2)$. However, in this case,  taking
$J=\{1,2\}$ and $J'=\{2,3\}$, one obtains $a_i=0$ for $i\in J\cap
J'=D$. Then $\mu=0$, which is false.

\med
Suppose that $G$ is twisted. We argue case-by-case.

\medskip
(i) $G={}^3D_4(q)$.
There are two $f$-orbits $J,J'$ on $D$, where 
 $J=\{1,3,4\}, J'=\{2\}$, and hence  $G_{J'}\cong SL(2,q)$
and $G_J\cong SL(2,q^3)$. So the assumption is not satisfied
unless $q=p$. This case occurs in item $(5)$. The claim on the
weights follows from Corollary \ref{fg8}.

\med
(ii)  $G={}^2A_{n}(q)$, where $n=2m-1>1$ is odd. Take $J=\{1\ld
m-1, m+1\ld 2m-1\}$. Then  $G_J\cong SL(m,q^2)$. By assumption,
$V_J$ is Steinberg. So $a_1=\cdots =a_{m-1}=a_{m+1}=\cdots
=a_{2m-1}=q-1$. Therefore, $a_m<q-1$. Take $J=\{m\}$. Then
$G_J\cong A_1(q)$, whence $q=p$ and $a_m=0$. The case $n=3$ is
recorded in (4). Let $n>3$. Take $J=\{m-1,m,m+1\}$. Then $G_J\cong
{}^2A_3(q)$, which is a contradiction.

\med
(iii) $G={}^2A_{n}(q)$, $n=2m$ even. Take $J=\{1\ld m-1, m+2\ld
2m\}$. Then $G_J\cong A_{m-1}(q^2)$. Then  $a_1=\cdots
=a_{m-1}=a_{m+2}=\cdots =a_{2m}=q-1$. Therefore,
$a_{m}a_{m+1}<(q-1)^2$. Then take $J=\{m,m+1\}$. Then  $G_J\cong
{}^2A_2(q)$, and hence $q=2$, $a_{m}=a_{m+1}=0$. The case $m=2$ is
recorded in item (6). Let $m>2$. Then for $J=\{m-1,m,m+1,m+2\}$ we
have $G_J={}^2A_{4}(2)$, which is a contradiction.

\med
(iv) $G={}^2D_{n}(q)$, $n>3$. Take  $J=\{n-1,n\}$. Then $G_J\cong
A_1(q^2)$, and the assumption implies $a_{n-1}=a_n=q-1$, and hence
$a_i<q-1$ for some $i<n-1$. 
Next take $J=\{1\ld n-2\}$, so $G_J=SL(n-1,q)$. This implies
$n=4$, $q=2$ and $a_1=a_2=0$. Finally, take $J=\{2,3,4\}$. Then
$G_J$ is of type ${}^2A_3(2)$, so this option is ruled out.

\med
(v) $G={}^2E_6(q)$. Then (using Bourbaki's ordering of the nodes)
the orbits of $f$ are $(1,6)$, $(3,5)$, $(2)$, $(4)$. Take
$J=\{1,3,4,5,6\}$. Then $G_J\cong {}^2A_5(q)$, which implies
$a_1=a_3=a_4=a_5=a_6=q-1$. So $a_2<q-1$. Next take $J=\{2,4\}$.
Then $G_J\cong A_2(q)$, and $V_J$ is not trivial. This is a
contradiction. \hfill $\Box$

\med
We close this section by two results which illustrate the use of
Propositions \ref{r55} and \ref{sd3}.

\begin{corol}\label{g77}
Let $G\in \{B_n(q),n>2,C_n(q),n>1,D_n(q), n>3,{}^{2}D_{n+1}(q),n>2
\}$. Let $J=\{1\ld n-1\}$.  Let $\Phi$ be a PIM with socle $V$ and
$S=\si_{G,G_J}(V)$. Then either $S$ is self-dual or $c_\Phi\geq
2c_\Psi$.
\end{corol}

Proof. Let $P_J$ be the standard parabolic subgroup of $G$
corresponding to $J$ and let $L:=L_J$ be the standard Levi
subgroup of $P_J$. Then $G_J\cong SL(n,q)$. Let  $\overline{W}_L$
be as in Proposition \ref{r55}. We show first that
$|\overline{W}_L|=2$, and a non-trivial element of
$\overline{W}_L$ induces the duality \au $h$ on $L$ (that is, if
$\phi\in \Irr L$ then $\phi^h$ is the dual of $\phi$). Indeed, we
can assume that $T_0\subset L$, and then $N_G(T_0)$ contains an
element $g$ acting on $T_0$ by sending every $t\in T_0 $ to
$t\up$. Let $h$ be the inner \au of $G$ induced by the
$g$-conjugation.  Then $S^h$ is the dual $L$-module of $S$. So the
result follows from Proposition \ref{r55}. (One can make this
clear by using an appropriate basis of the underlying space $V$ of
$G$, and by considering, for the group $G_1$ of isometries of $V$,
a matrix embedding $GL(n,q)\ra G_1 $ sending every $x\in GL(n,q)$
to $\diag(x,{}^{Tr}x\up)$, where ${}^{Tr}x$ denotes the transpose
of $x$.)

\med
A similar result holds for the groups
${}^{2}A_{2n+2}(q),{}^{2}A_{2n+1}(q)$, where $P$ has to be chosen
so that $L$ contains $SL(n,q^2)$. In this case $h$ is the duality
automorphism following the Galois automorphism. 

\begin{propo}\label{bn1}
Let $k=0,1$ and $G={}^2A_{2n-1+k}
\cong SU(2n+k,q)$. Let $J=\{1\ld n-1,n+k+1\ld 2n-1+k\}$ so
$G_J\cong SL(n,q^2)$. Let $\mu=a_1\om_1+\cdots
+a_{2n-1+k}\om_{2n-1+k}\in \Delta(G)$, let $V_\mu$ be an \ir  $\GC$-module of
highest weight $\mu$, $V=V_\mu|_G$ and $V_J=\pi_{G,G_J}$.

\med
$(1)$ The highest weight of $V_J$ is
$\mu'=(a_1+qa_{2n-1+k})\om_1'+\cdots+(a_{n-1}+qa_{n+k+1})\om_{n-1}'
$, where $\om'_1\ld \om_{n-1}'$ are the fundamental weight of
$SL(n,F)$.

\med
$(2)$ Let $\Phi$ be the PIM with socle $V$, and $\Psi$ be the PIM
for $G_J$ with socle $V_J$. Then either $c_{\Phi}\geq 2c_{\Psi}$,
or $a_i=a_{n+i+k}$ for $i=1\ld n-1$. 
\end{propo}

Proof. (1) follows from Proposition \ref{sd3}. (2) Let $1\neq w\in
\overline{W}_L$. The  \au induced by $w$ on $G_J$ sends $g\in
SL(n,q^2)$ to ${}^{Tr}g^{-\gamma}$, where $\gamma$ is the Galois
\au of $F_{q^2}/F_q$ and ${}^{Tr}g$ is the transpose of $g$. Then
the mapping $v\ra w(g)v$ $(v\in V_J)$ yields an \ir $FG_J$-module
 $V'_{J} $ of highest weight
$\mu':=(qa_{n-1}+a_{n+k+1})\om_1'+\cdots+(qa_1+a_{2n-1+k})\om_{n-1}'
$. If $V_{J} $ is not isomorphic to $V'_{J}$ then $c_\Phi\geq
2c_\Psi$ by Proposition \ref{r55}. If $V_J\cong V'_{J} $ then
$a_1+qa_{2n+k-1}\equiv (qa_{n-1}+a_{n+k+1})({\rm mod }\,q^2)\, \ld
\,$ $a_{n-1}+qa_{n+k+1}\equiv (qa_1+a_{2n+k-1})({\rm mod }\,q^2)$.
As $a_i<q$ for $i=1\ld 2n+k$, it follows that
$a_i=a_{n+k+i}$ for $i=1\ld n-1$.  \hfill $\Box$

\section{Harish-Chanra induction and a lower bound for the PIM dimensions}

The  classical result by Brauer and Nesbitt \cite[Theorem 8]{BN}
(see also \cite[Ch.IV, Lemma 4.15]{Fe}) is probably not  strong
enough to be used  for  studying PIM dimension bounds for
Chevalley groups.

\med
The \f lemma is one of the standard results on \reps of groups
with $BN$-pair, see  \cite[pp.683 - 684]{CR2} or \cite[\S 14,
Theorem 48]{St}:

\begin{lemma}\label{bb5}
Let $G$ be a Chevalley group
viewed as a group with BN-pair, $B$ a Borel subgroup and $W_0$ the
Weyl group. Then the induced character $1_B^G$ is a  sum of \ir
characters $\chi_\lam$ of $G$ labeled by $\lam \in \Irr W_0$, and
the multiplicity  of $\chi_\lam$ in $1_B^G$ is equal to $\lam
(1)$. In other words, $1_B^G=\sum_{\lam\in\Irr
W_0}\lam(1)\chi_\lam$.
\end{lemma}

Recall that $G$ is called  a Chevalley group if $ G =\GC^F$, where
$\GC$ is a simple, simply connected algebraic group.

\begin{propo}\label{r01a}
Let $G$ be a  Chevalley group,
$B$ a Borel subgroup of $G$, $U$ the \syl of $B$ and $T_0$  a
maximal torus of $B$. Let  $W_0$ be the Weyl group of $G$ as a
group with $BN$-pair and let $d$ be the minimum degree of a
non-linear
character of $W_0$. Let $\chi$ be a 
 character vanished at all $1\neq u\in U$  such that $(\chi,1_G)=(\chi,St)=0$.
 Suppose that
the derived subgroup of $W_0$ has index $2$ and, for every
non-trivial  character $\beta\in \Irr T_0$, either $|W_0\beta|\geq
d$ or $(\chi,\beta_B^G)=0 $, where $\beta_B$ is the inflation of
$\beta$ to $B$. Then  $\chi(1)\geq d\cdot |G|_p$.

\end{propo}

 Proof. By Corollary \ref{gg1}, $\chi(1)=(\chi,1_U^G)\cdot |G|_p$.
 By Lemma \ref{m12}, $(\chi,1_U^G)= \sum_{\beta_B}|W_0\beta|(\chi,
\beta_B^G)$. Suppose $\chi(1)<d\cdot |G|_p$. Then, by assumption,
$(\chi, 1_U^G)=(\chi,1_B^G)$.

Let $\chi_\lam$ be the \ir constituent of $1_B^G$ corresponding to
$\lam\in\Irr W_0$, see Lemma \ref{bb5}. Let $m_\lam$ be the
multiplicity of $\chi_\lam$ in $\chi$. By Lemma \ref{bb5},
$$(\chi,1_B^G)=\sum_{\lam\in { \Irr W}}m_\lam \cdot \lam (1)<d .$$
\itf $m_\lam =0$ if $\lam(1)>1$. So $(\chi,1_B^G)=\sum_{\lam\in
\Irr W:\lam(1)=1}m_\lam$. It is well known \cite[Theorem
67.10]{CR2} that $1_G$ and $St$ occur in $1_B^G$ with \mult 1. As
$W_0$ has exactly two one-dimensional representations, $\chi_\lam$
is either $St $ or $1_G$. This is a contradiction, as neither $St$
nor $1_G$ is a constituent of $\chi$. \hfill $\Box$

\med
Remarks. (1) The values of $d$ are given by Table 2. (2)  The
derived subgroup of $W_0$ has index 2 \ii $G\in \{A_n(q),D_n(q),
E_6(q)$, $E_7(q)$, $E_8(q)$, ${}^2B_2(q)$, ${}^2G_2(q)\}$. (The
last two groups are of BN-pair rank 1, and hence $|W_0|=2$.)

\begin{propo}\label{r01}
Let $G$ be a  Chevalley group, $B$ a Borel subgroup of $G$, and
$T_0$  a maximal torus of $B$. Let  $W_0$ be the Weyl group of $G$
as a group with $BN$-pair and let $d$ be the minimum degree of a
non-linear character of $W_0$. Let $\Phi\neq St$ be a PIM for $G$
with character $\chi$. Suppose that the derived subgroup of $W_0$
has index $2$ and $|W_0\beta|\geq d$ for every non-trivial
character $\beta\in \Irr T_0$,
where $\beta_B$ is the inflation of $\beta$ to $B$. Then
$c_\Phi\geq d$, unless $\Phi=\Phi_1$ and $G=SL(2,p)$, $SL(3,2)$ or
${}^2G_2(3)$.
\end{propo}

 Proof. We show that the hypothesis of Proposition \ref{r01a} holds.
 Indeed, $\chi (u)=0$ for every $1\neq u\in U$, as $\chi$ is the character
 of a projective module. As $St$ is a character
of a PIM, $St$ is not a constituent of  $\chi$. Suppose that
$(\chi,1_G)>0$.
 Then, by orthogonality relations \cite[Ch. IV, Lemma 3.3]{Fe},
 we have $\Phi= \Phi_1$ and
$(\chi,1_G)=1$. So $\dim\Phi_1=|G|_p$. This contradicts a result
of Malle and Weigel \cite{MW}, unless $G= SL(2,p) $, $SL(3,2)$  or
${}^2G_2(3)$. With these exceptions, the result now follows from
Proposition \ref{r01a}. \hfill $\Box$

\begin{theo}\label{th1}
 Let $G\in \{SL(n,q)$, $n>4$,  $ Spin^{+}(2n,q)$, $q$ even, $n>3$,
$E_6(q),E_7(q),E_8(q)\}$. Let $\Phi\neq St$ be a PIM with socle
$S$.
\med

$(1)$ If $\soc S|_B\neq 1_B$ then $ c_\Phi\geq m$, where
$m=n,2n,27,28,120$, respectively.

\med
$(2)$ Suppose that $\soc S|_B= 1_B$. Then $c_\Phi\geq d$, where
 $d=n-1,n-1,6,7,8$, respectively.
 \end{theo}

Proof. Let $W$ be the Weyl group of $G$.

 (1) By Proposition \ref{r2}, $c_\Phi\geq |W\beta|$. The
lower bound $m$ for $|W\beta|$ is provided in  Propositions
\ref{r4} and \ref{dn5}. In particular,  $|W\beta|\geq m$ unless
$C_W(\beta)=W$. This implies $\beta=1_{T_0}$ by  Lemma \ref{cm2}.

\med
(2) follows from  Proposition \ref{r01}. Indeed, let $d$ be  the
minimum dimension of a non-linear \irr of $W$; by  Table $2$, $d$
is as in the statement (2). As $m> d$, it follows from (1) that
$|W\beta|> d$ for every $1_B\neq \beta\in \Irr T_0$. For the
groups $G$ in the statement the derived subgroup of $W$ is well
known to be of index 2. So the hypothesis of Proposition \ref{r01}
holds, and so does the conclusion. \hfill $\Box$

\begin{corol}\label{or4}
Let $G=
 {\rm Spin}^+\,(2n,q),$
 $ n>4$, and
let $\Phi\neq St$ be a PIM with socle $V=V_\mu|_G$, where $V_\mu$
is a $\GC$-module with highest weight $\mu =a_1\om_1+\cdots
+a_n\om_{n}$. Then $c_\Phi\geq n-1$.
\end{corol}

Proof.
 Let  $J=\{1\ld n-1 \}$ or $  \{1\ld n-2,n \}$.
 Then
$G_J\cong SL(n,q)$ and the highest weight of $V_j$ is
$\nu:=a_1\om_1+\cdots +a_{n-1}\om_{n-1}$ or $a_1\om_1+\cdots
+a_{n-2}\om_{n-2}+a_{n}\om_{n}$, respectively. So $V_J$ is not the
Steinberg $FG_J$-module at least for one of these two cases.
Furthermore, by Theorem \ref{th1}, if $V_J\neq St$ then
$c_{\Psi}\geq n-1$, where $\Psi$ is a PIM for $G_J$ with socle
$V_J$. By Lemma \ref{mp1}, $c_{\Phi} \geq c_{\Psi}$, and the
result follows. \hfill $\Box$

\med
Remark. If $G= SL(n,q)$ or Spin${}^+(2n,q)$ then it follows from
Theorem \ref{th1} and Corollary \ref{or4} that $c_\Phi\ra \infty$
as $n\ra \infty$, provided $\Phi\neq St$.
  If $G\in \{Sp(2n,q), Spin(2n+1,q) $, Spin${}^-(2n+2,q)\}$ then a similar argument
gives that $c_\Phi\geq n-1$ unless $a_1=\cdots =a_{n-1}=q-1$
(respectively, $c_\Phi\geq n-2$ unless $a_1=\cdots =a_{n-2}=q-1$).
However, this does not lead to the same conclusion as above. A
similar difficulty arises for the groups  $SU(2n,q)$ and
$SU(2n+1,q)$.

\begin{propo}\label{tp1}
Let $G=G(q)$ be a non-twisted Chevalley group, and let $V_\mu$ be
a $\GC$-module with \hw $\mu=\sum a_i\om_i \in \Delta(G)$. Let
$\Phi$ be a PIM of $G$ with socle $V_\mu|_G$. Let $J$ be a subset
of  nodes on the Dynkin diagram of $\GC$, not adjacent to each
other. \st $a_i\neq 0,q-1$ for $i\in J$.
 Then $c_\Phi\geq
2^{|J|} $.
\end{propo}

Proof. As the nodes  of $J$ are not adjacent, it follows that
$\GC_J$ is the direct product of $|J|$ copies of $G_i\cong
SL(2,F)$ for $i\in G$. Furthermore, the Smith correspondent
$\sigma_{G,G_J}$ of $V_\mu$ is the tensor product of \ir
$FG_i$-modules  with highest weight $a_i\om_1$.  Let
$\Psi=\pi_{G,G_J}(\Phi)$ be the parabolic descent of $\Phi$ to
$G_J$. Then $c_{\Psi}\geq 2^{|J|} $ by Lemmas \ref{sr1},
\ref{dp8}, and $c_{\Phi} \geq c_{\Psi}$ by parabolic descent. \hfill $\Box$

\med

The proof of Proposition \ref{tp1}  illustrates the fact that the
parabolic descent from  non-twisted groups $G(q)$ for $q=p$ to
minimal parabolic subgroups (distinct from the Borel subgroup)
does not work. Indeed, in this case $G_J\cong A_1(p)$;  if all
coefficients $a_i$ are equal to 0 or $p-1$ then $c_\Psi=1$. One
observes that it is possible to run the parabolic descent to
subgroups $G_J\cong SL(3,p)$, obtaining tensor product of \ir
\reps with highest weight $(p-1)\om_1$ or $(p-1)\om_2$, and then
use known results for PIM's with such socles.

\section{The Ree groups}

In this section we prove Theorem \ref{wt1} for Ree groups
$G={}^2G_2(q)$, $q^2=3^{2f+1}>3$.  Note that
$|G|_p=q^{12}=3^{6(2f+1)}$.

Note that if $f=0$ then $G\cong {\rm Aut}\,\SL(2,8)$. For this
group the decomposition numbers are known, and $c_\phi=1$ \ii
$\Phi=\Phi_1$.

\begin{lemma}\label{re1} 
Let $\chi$ be a  character of degree $q^6$ vanishing at all
unipotent elements of $G$. Suppose that $(\chi,1_G)=0$. Then
$\chi=St$.
\end{lemma}

Proof. Suppose the contrary.
 Then $(\chi, 1_U^G)=1$ (Corollary \ref{gg1}),  so there is exactly one
\ir \ccc $\tau$ of $\chi$ common with $1_U^G$. Recall that
$1_U^G=\sum_{\beta\in \Irr T_0 }\beta_B^G$, where $\beta_B$ is the
inflation of $\beta$ to $B$. By Proposition \ref{r2b}, $\tau\in
\beta_B^G $ implies that $\beta$ is $W_0$-invariant. (Or
straightforwardly, the degree of every $\beta_B^G$ is
$|G:B|=|G|_p+1$, so we can ignore the $\beta$'s with  $\beta_B^G$
irreducible. Observe that $\beta_B^G$ is reducible \ii
$C_{W_0}(\beta)=W_0$ \cite[Theorem 47]{St}.) As $T_0$ is cyclic of
order $q^2-1$ and $|W_0|=2$, the non-identity element of $W_0$
acts on $T_0$ as $t\ra t\up$ ($t\in T_0$). So $C_W(\beta)\neq 1$
implies $\beta^2=1$. In this case
$(\beta_B^G,\beta_B^G)=|C_W(\beta)|=2$ (\cite[Ex.(a) after Theorem
47]{St}), so there are two \ir constituents. If $\beta=1_{T_0}$
then $\beta_B^G=1_G+St$. If $\beta^2=1_{T_0}\neq \beta$ then the
character table in \cite{Wa} leaves us with exactly two
possibilities for $\tau(1)$, which are $d_1=q^4-q^2+1$ and
$d_2=q^6+1-d_1=q^6-q^4+q^2$.

 The constituents of $\chi$ that are not in $1_U^G$ 
 are cuspidal. (Indeed,  every proper parabolic subgroup of $G$ is a Borel
subgroup, so every non-cuspidal \ir character  has a 1-dimensional
constituent under restriction to $U$. So they are in $1_U^G$.)
Malle  and Weigel \cite[p. 327]{MW} recorded the cuspidal
character degrees that are less than $|G|_p$, and each degree
 is greater than $q^6-d_2$. Therefore, $\tau(1)=d_1$. By Lemma
\ref{gg1}, there is exactly one regular character $\gamma$
occurring  as a constituent of $\chi$, and $\gamma$ is not a
unipotent character (as $St$ is the only unipotent regular
character of $G$). The two  characters listed in \cite[p. 327]{MW}
are unipotent (see \cite[p. 463]{Ca} for the degrees of unipotent
characters of $G$). The remaining regular character degrees are

\med
$d_3=(q^2-1)(q^4-q^2+1)=q^6-2q^4+2q^2-1$, and

\med
$d_4=(q^4-1)(q^2-q\sqrt{3}+1)=q^6-q^5\sqrt{3}+q^4-q^2+q\sqrt{3}-1$.

\med
 Thus,  $\gamma(1)=d_3$ or $d_4$. As $(\chi, 1_U^G)=1$, we have $(\tau,
 \chi)=1$.  Note that $q^6-d_1-d_3=q^4-q^2$.
    Other characters which may occur in the decomposition of $\chi$
  are  of degree $d_5=q(q^4-1)/\sqrt{3}$,  $d_6=q(q^2-1)(q^2-
  q\sqrt{3}+1)/2\sqrt{3}$ or $d_7=q(q^2-1)(q^2+
  q\sqrt{3}+1)/2\sqrt{3}$, see \cite{Wa}. Each of these degrees
  is greater than $q^4-q^2$, so $\gamma(1)=d_4$.
So

\med
$$q^6=(q^4-q^2+1)+ (q^6-q^5\sqrt{3}+q^4-q^2+q\sqrt{3}-1)+$$
\begin{equation}\label{e1}
+\frac{aq(q^4-1)}{\sqrt{3}}+ \frac{b q(q^2-1)(q^2+ q\sqrt{3}+1)}{2\sqrt{3}}
+\frac{cq(q^2-1)(q^2-  q\sqrt{3}+1)}{2\sqrt{3}},\end{equation}

\med
\noindent where $a,b,c\geq 0$ are integers, and $a+b+c>0$.
\med

Cancelling the equal terms and  dropping
the common multiple $q^2-1$, we get

 \med
\begin{equation}\label{e3}
a(q^2+1)  + (q^2+1)\frac{b+c}{2}-3(q^2+1)=q\sqrt{3}\cdot (\frac{c-b}{2}-2).\end{equation}

\med
Suppose first that $a>2$. Then $(q^2+1)\frac{b+c}{2}\leq
q\sqrt{3}\cdot (\frac{c-b}{2}-2) $, and hence $e:=c-b>0$. Then we
have $(q^2+1)b+e\cdot (\frac{q^2+1}{2}-\frac{q\sqrt{3}}{2})\leq
-2q\sqrt{3}$. One easily check that $q^2+1-q\sqrt{3}>0$, which is
a contradiction. So $a\leq 2$.

\med
Let   $X,Y,m$ be as in the character table of $G$ in \cite{Wa}, so
 $m=q/\sqrt{3}$, $|X|= 3$ and $|Y|=9$. Then one observes from
\cite{Wa} that the \ir characters of the same degree have the same
value at $X$, and also at $Y$. Therefore,
$\chi(g)=\tau(g)+\gamma(g)+a\xi_5(g)+b\xi_6(g)+c\xi_7(g)=0$, where
$g\in\{X,Y\}$ and $\xi_j(g)$ is the value of any character of
degree $d_j$ ($j=5,6,7$). By  \cite{Wa}, $\tau (Y)=-\gamma (Y)=1$,
$\xi_5(Y)=-m$, $\xi_6(Y)=\xi_7(Y)=m$. So
$a\xi_5(g)+b\xi_6(g)+c\xi_7(g)=-am+(b+c)m=0$, whence $a=b+c$. So
(\ref{e3}) simplifies to

\begin{equation}\label{e2}
3(q^2+1)((b+c-2)=q\sqrt{3}\cdot (c-b-4).\end{equation}

\med
Recall that $a=b+c\leq 2, $ and hence $b+c\neq 0$. If $b+c=2$ then
$c=b+4$, which is a contradiction. If $b+c=1$ then the left hand
side in $(\ref{e2})$ is not divisible by 9, and hence
$q\sqrt{3}=3$. This is not the case. \hfill $\Box$

\begin{propo}\label{re2} Let
$G={}^2G_2(3^{2f+1})$, $q^2=3^{2f+1}$. Let $\Phi\neq St$ be a PIM,
and let $\chi$ be the character of $\Phi$. Then $\chi( 1)>q^6$.
\end{propo}

Proof. If  $(\chi,1_G)=0$, the result follows from Lemma
\ref{re1}. Otherwise $\Phi=\Phi_1$, and the result for this case
follows from \cite{MW}. \hfill $\Box$

\med
Proposition \ref{re2} together with certain known results implies
Theorem \ref{wt1} for groups of BN-pair rank 1. These are $A_1(q),
{}^2A_2(q)$, ${}^2B_2(q)$, ${}^2G_2(q).$ The groups ${}^2G_2(q)$
have been treated above. Below we quote known results on minimal
PIM dimensions for the remaining groups of BN-pair rank 1.

\begin{lemma}\label{sr1}
{\rm \cite{Sr}} Let
$G=SL(2,q)$, $q=p^k$, and let $\Phi\neq St$ be a PIM of $G$ with
socle $V_{m\om_1}|_G$, where $V_{m\om_1}$ is an \ir
$SL(2,F)$-module with highest weight $m\om_1$ for $0<m<q-1$. Let
$m=m_0+m_1p+ \cdots+ m_{k-1}p^{k-1}$ be the $p$-adic expansion of
$m$, and $r$ the number of digits $m_i$ distinct from $p-1$. Then
$c_\Phi=2^r$. In addition, $c_{\Phi_1}=2^k-1$.  In particular, if
$q=p$ then $c_\Phi=2$ unless $\Phi=\Phi_1$ with $c_{\Phi_1}=1$.
\end{lemma}

\begin{lemma}\label{sr1a}
{\rm \cite{CF1}}
Let $G={}^2B_2(q)$, $q^2>2$, and let $\Phi\neq St$ be a PIM of
$G$. Then $c_\Phi\geq 4$. \end{lemma}

\begin{lemma}\label{su0} Let $G=SU(3,p)$, $p>2$, and let $\Phi\neq St$ be a PIM of $G$.
 Then $ c_\Phi \geq 3 $.
\end{lemma}

Proof. If $p\equiv 1 \,({\rm mod}\, 3)$ then the PIM dimensions
are listed in   \cite[Table 1]{Hu90}.  The same holds if $p=3$,
see the GAP library. Let $p\equiv 2  \,({\rm mod}\, 3)$. Then
 the formulas in \cite[p. 10]{Hu81}, lead to the same conclusion.

\begin{propo}\label{ra1} Theorem
$\ref{wt1}$ is true for groups of BN-pair rank $1$.
\end{propo}

Proof. The result follows from Proposition \ref{re2} and Lemmas
\ref{sr1}, \ref{sr1a} and \ref{su0}. \hfill $\Box$

\section{Groups $SU(4,p)$ and ${}^3D_4(p)$}

In this section we prove Theorem \ref{wt1} for the above groups.

 \subsection{Some general observations}\label{74}

Prior dealing with the cases where $G=SU(4,p)$ and ${}^3D_4(p)$,
we make some general comments which facilitate computations. These
are valid for $q$ in place of $p$, so we do not assume $q$ to be
prime until Section 8.2.

\med
Let ${\mathbf H}$ be a connected reductive algebraic group and
$H={\mathbf H}^{Fr}$. The Deligne-Lusztig theory partitions \ir
characters of $H$ to series  usually denoted by ${\mathcal E}_s$,
where $s$ runs over representatives of the semisimple conjugacy
classes of the dual group $H^*$. (See \cite[p. 136]{DM}, where our
${\mathcal E}_s$ are denoted by ${\mathcal E}(\GC^F,
(s)_{\GC^{*F^*}}).)$ The duality also yields a bijection $T\ra
T^*$ between maximal tori in $H$ and $H^*$  such that $T^*$ can be
viewed as $\Irr T$, the group of linear characters of $T$. If
$H=U(n,q)$ or ${}^3D_4(q)$ then $H^*\cong H$.

\begin{lemma}\label{dd4}
Let $H$ be a finite reductive group, $B$  a Borel subgroup of $H$
and $T_0\subset B$  a maximal torus. Let $B=T_0U$, where $U$ is
the unipotent radical of $B$. Let $\beta\in\Irr T_0$, and
$\beta_B$ its inflation to $B$. Let $s\in T_0^*\subset H^*$ be the
element corresponding to $\beta$ under the duality mapping $\Irr
T_0\ra T_0^*$.

\medskip
$(1)$ ${\mathcal E}_s$ contains the set of all \ir constituents of
$\beta_B^H$. In addition, if an \ir character $\chi$ of $H$ is a
constituent of $1_U^H$ then $\chi\in{\mathcal E}_s$ for some $s\in
T_0^*$.

\medskip
$(2)$  $\beta_B^H$ contains a  regular character and a semisimple
character of ${\mathcal E}_s$.

\medskip
$(3)$ Suppose that  ${\mathcal E}_1$ coincides with the set of all
\ir constituents of $1_B^G$ and $s\in Z(H^*)$. Then ${\mathcal
E}_s$ coincides with the set of all \ir constituents of
$\beta_B^G$.

\end{lemma}

Proof. (1) By  \cite[Proposition 7.2.4]{Ca},
$\beta_B^H=R_{T_0,\beta} $, where $R_{T_0,\beta} $ is a
generalized character of $H$ defined by Deligne and Lusztig
(called usually a Deligne-Lusztig character). By definition
\cite[p. 136]{DM}, ${\mathcal E}_s$ contains the set  $\{\chi\in
\Irr H:(\chi,R_{T_0,\beta} )\neq 0\}$, which coincides with
$\{\chi\in \Irr
H:(\chi,\beta_B^H )\neq 0\}$. 
This implies the first claim. (Note that this statement is a
special case of \cite[Proposition 4.1]{Ko}.) The additional claim
follows as $1_U^H=\sum_{\beta\in\Irr T} \beta_B^H$.

\medskip
(2) This is a special case of   \cite[Proposition 5.1 and
Corollary 5.5]{Ko}.

\medskip
(3) By \cite[Proposition 13.30]{DM}, there exists a
one-dimensional character $\hat s$ of $H$ such that the characters
of ${\mathcal E}_s$ are obtained from those of ${\mathcal E}_1$ by
tensoring with $\hat s$, and $\hat s|_{T_0}=\beta$. It follows
that ${\mathcal E}_s$ is the set of all \ir constituents of
$1_B^H\otimes \hat s$. As $\hat s|_B=\beta_B$, we have
$1_B^H\otimes \hat s=\beta_B^H$, as required. \hfill $\Box$

\medskip
Remark. (3) can be stated by saying that if ${\mathcal E}_1$
coincides with the Harish-Chandra series of $1_B^H$ and $s\in
Z(H^*)$ then ${\mathcal E}_s$ coincides with the Harish-Chandra
series of $\beta_B^H$.

\begin{lemma}\label{gg4} $(1)$ Let $s\in G^*$ be a semisimple element
and ${\mathcal E}_s$ the corresponding Lusztig class. Let $\Gamma$
be a Gelfand-Graev character of $G$. Then ${\mathcal E}_s$
contains exactly one character common with $\Gamma$.

$(2)$ Suppose that $s\in T^*_0=\Irr T_0\cong T_0$. Let
$\beta\in\Irr T_0$  correspond to $s$. Then every regular
character of ${\mathcal E}_s$ is a constituent of  $\beta_B^G$.
\end{lemma}

Proof. (1) is well known if $\GC$ has connected center. In general
this follows from results in \cite[Section 14]{DM}. Indeed, the
character  $\chi_{(s)}$ introduced in \cite[14.40]{DM} is a sum of
regular characters \cite[14.46]{DM}. By definition of
$\chi_{(s)}$, its \ir constituents belong to ${\mathcal E}_s$. In
addition, every regular character of  $G$ is a constituent of
$\chi_{(s)}$ \cite[14.46]{DM}. Therefore, every regular character
of ${\mathcal E}_s$ is a constituent of $\chi_{(s)}$. As $(
\chi_{(s)}, \Gamma)=1$ \cite[14.44]{DM}, the claim follows.

(2) By \cite[7.4.4]{Ca}, $\beta_B^G=R_{\TC_0,\beta}$, where
$R_{\TC_0,\beta}$ is the Deligne-Lusztig character. Note that
$(R_{\TC_0,\beta}, \Gamma)=1$ (indeed, it follows from the
reasoning in the proof of Lemma 14.15 in \cite[14.44]{DM} that
$(R_{\TC_0,\beta}, \Gamma)=\pm 1$; as $\beta_B^G$ is a character,
$(\beta_B^G,\Gamma)\geq 0$).  (This also follows from
\cite[Proposition 2.1]{Ko2}.) So the claim follows from (1).\hfill
$\Box$


\begin{lemma}\label{nf1}
 Let $G=SU(4,q)$, resp., ${}^3D_4(q)$, and let $J=\{1,3\}$, resp., $\{1,3,4\}$
be a set of nodes on the  Dynkin diagram of the respective algebraic group
 $\GC$.
Let $\mu=(q-1)\sum_{j\in J}\om_j$,  $V_\mu$ an \ir $\GC$-module of
 highest weight $\mu$, and $V=V_\mu |_G$. Then $S:=\si _{G,L}(V)$
 is of defect $0$. Let $\Psi$ be a PIM with socle $S$ and let
  $\lam$ be the character of the lift of $S$. Then $(\lam^{\# G},\chi)>0$.
  In addition, if $c_\Phi=1$ then $(\lam^{\# G},\tau)=1$.
\end{lemma}

Proof. The statement about the defect of $S$ and the inequality
$(\lam^{\# G},\chi)>0$ is proved in Proposition \ref{yg8}(3).
Moreover, it is shown there that $\lam$ is a constituent of
$1_{B\cap L}^L$. Therefore, by transitivity of Harish-Chandra
induction \cite[Proposition 70.6(iii)]{CR2}, $1_B^G=\lam^{\#
G}+\lam'$, where $\lam'$ is some character of $G$. So $\tau$ is a
common constituent of $\lam^{\# G}$ and $1_U^G$; as
$(\tau,1_U^G)=1$, it follows that $(\tau, \lam^{\# G})=1$. \hfill
$\Box$

\begin{lemma}\label{u44}
Let $H=U(4,q)$, $B$  a Borel subgroup of $H$ and $T_0\subset B$ a
maximal torus. Let $B=T_0U$, where $U$ is the unipotent radical of
$B$. Let $\beta\in\Irr T_0$, and $\beta_B$ its inflation to $B$.
Let $s\in T_0^*\subset H^*$ be the element corresponding to
$\beta$ under the duality mapping $\Irr T_0\ra T_0^*$.

\medskip
$(1)$ ${\mathcal E}_1$ coincides with the set of all \ir
constituents of $1_B^H$.

\medskip
$(2)$ ${\mathcal E}_s$ coincides with the set of all \ir
constituents of $\beta_B^H$.
\end{lemma}

Proof. (1) Note that $H$ has no cuspidal unipotent \ir character,
see \cite[p. 457]{Ca}. \itf every unipotent character of $H$ is a
constituent of $1_B^G$. This is equivalent to the statement.

\medskip
(2) Observe  that $H\cong H^*$. In notation of \cite{No},
$A_1,A_9,B_1,C_1,C_3$ are the only conjugacy classes that meet
$T^*_0$. If $s\in A_1$ then $s\in Z(H^*)$, and the claim follows
Lemma \ref{dd4}(2). If $s$ belongs to $B_1,C_1,$ or $C_3$ then
  ${\mathcal E}_s$ contain only regular
and semisimple characters, so the claim follows from (2). Let
$s\in A_9$. By \cite[13.23]{DM}, $|{\mathcal E}_s|$ is equal to
the number of unipotent characters of $C_{H^*}(s)$. This group is
isomorphic to $U (2,q)\times U(2,q)$, so the number of unipotent
characters equals 4. By Lemma \ref{dd4}(1), the \ir constituents
of $\beta_B^H$ are contained in ${\mathcal E}_s$. We show that
$\beta_B^H$ has 4 distinct \ir constituents.

Indeed, $(\beta_B^H,\beta_B^H)=|W_s|$, where $W_s$ is the
stabilizer of $s$ in  $W_0$ (see Section 4 for the definition of
$W_0$). It is easy to observe that $W_s$ is an elementary abelian
group of order 4. By (2), the regular and semisimple characters
occur as  constituents of  $\beta_B^H$. As
$(\beta_B^H,\beta_B^H)=4$ is the sum of squares of the
multiplicities of the \ir constituents of $\beta_B^H$, it follows
that  there are exactly four distinct constituents. This equals
$|{\mathcal E}_s| $, the result follows.\hfill $\Box$

\med
Let $\Phi$ be  a PIM  for $G$ with character $\chi$.
 Suppose $c_\Phi =1$. Then, by Corollary \ref{gg1},
   $(\chi, 1_U^G)=(\chi, \Gamma)=1$, which means that $\chi$ has exactly
   one \ir constituent common with $1_U^G$, and  exactly one \ir
   constituent common with every Gelfand-Graev character $\Gamma$.
   Denote them as $\tau$ and $\gamma$, respectively. By Proposition \ref{yg8},
$(\chi, 1_B^G)>0$; as $(\chi, 1_U^G)=1 $, it follows that $\tau$
is constituent of $1_B^G$. (This also follows from Proposition
\ref{r2b}.) Recall that the \ir constituents of $\Gamma$ are
called regular characters in
  the Deligne-Lusztig theory.  
  In general
there are several Gelfand-Graev characters, however, if $G=U(n,q)$
 or ${}^3D_4(q)$ then
$G$ has a single  Gelfand-Graev character (see \cite[14.29]{DM}).
 Note that by Lemma \ref{de4}, it suffices to prove that $c_\Phi>1$
for the group $U(4,q)$ in place of $SU(4,q)$.

\med
In order to prove Theorem \ref{wt1} for the   groups $G=SU(4,p)$
and ${}^3D_4(p)$ we first determine $\tau(1)$ and $\gamma(1)$, and
observe that $\tau\neq \gamma$ by Lemma \ref{dd4}. Next we express
$\chi=\tau+\gamma+\sum \nu_i$, where $\nu_i$ runs over the
characters that are neither regular nor in $1_U^G$. In particular,
  $\nu_i(1)\leq \chi(1)-\tau(1)-\gamma(1)$.
As the character table of $G$ is available in \cite{No} for
$U(4,q)$ and in \cite{DM} for ${}^3D_4(q)$, we obtain a
contradiction by inspecting all the possibilities.

The reasoning below does not use much from modular \rep theory. In
fact, we prove the following. Let $G$ be either $SU(n,q)$ or
${}^3D_4(p)$, and let $\chi$ be an ordinary character vanishing at
all non-semisimple elements of $G$. Let $L$ be a Levi subgroup of
$G$ whose semisimple part is isomorphic to $SL(2,q^a)$, where
$a=2$ for the former group and $a=3$ for the latter one. Suppose
that $\overline{\chi}_L$ is of defect zero. Then $\chi(1)>|G|_p$.
(Of course one has to replace $p$ by $q$, and use Lemma \ref{b1}
in place of Corollary \ref{gg1}.)

\subsection{The unitary groups $G=SU(4,p)$}


\med
Note that $|G|=p^6(p^4-1)(p^3+1)(p^2-1)$. Let $\Phi$ be a PIM with
socle  $V$ and character $\chi$. Arguing by contradiction, we
suppose that $\dim\Phi=|G|_p$, that is, $c_\Phi=1$.

 Let $J=\{1,3\}$ be a
set of nodes at the Dynkin diagram of $\GC=SL(4,F)$. Then
$G_J\cong SL(2,p^2)$, see Example prior Lemma \ref{st0}. Let
$P=P_J$  be a  parabolic subgroup of $G$ corresponding to $J$ and
$L=L_J$ its Levi subgroup. Let $S=\si_{G,P_J}(V)$ be the
Smith-Dipper correspondent of $V$. By Corollary \ref{fg8},
$S|_{G_J}$ is the Steinberg module. So $S$ is an $FL$-module of
defect 0, and hence lifts to characteristic 0. Let $\lam$ be the
character of the lift so $\lam(1)=p^2$. By Lemma  \ref{nf1},
$\tau$ is a constituent of $\lam^{\# G}$.

\med
In order to determine $\tau(1)$ we first decompose $\lam^{\# G}$
as a sum of \ir constituents.

\begin{lemma}\label{u4st}
$(1)$ $\lam^{\# G}=St+\si+\si'$, where $\si,\si'$ are \ir
 characters of $G$ such that
$\si (1)=p^2(p^2+1)$ and $\si'(1)=p^3(p^2-p+1)$.

\med
$(2)$ Let $\tau$ be a common \ir constituent for $\chi$ and  $
1_U^G$.  Then $\tau=\si'$, where $\si'$ is as in $(1)$. In
particular, $\tau(1)=p^3(p^2-p+1)$.
 \end{lemma}

Proof. (1) The degrees of \ir constituents of $1_B^G$ are given in
\cite[Proposition 7.22]{CKS}. As $\lam^{\#
G}(1)=p^2|G:P_J|=p^2(p+1)(p^3+1)$, this implies (1).

\med
 (2) By Corollary \ref{gg1}, $(\chi,1_U^G)=1$; as mentioned
  in \cite[Proposition 7.22]{CKS},  $\si$ occurs
  in $1_B^G$ with \mult
  2, so  $(\chi, \si)=0$, and hence    $\tau=\si'$, as stated.
 \hfill $\Box$

\med

 We need to write down the degrees of the regular characters of
$H:=U(4,p)$.  By the Deligne-Lusztig theory,   the regular
characters of $H$ are in bijection with  semisimple conjugacy
classes in the dual group $H^*\cong H=U(4,p)$. So we write
$\rho_s$ for the regular characters of $H$ corresponding to $s$, a
representative of a semisimple conjugacy classes in $H$.
Furthermore, $\rho_s(1)=|C_H(s)|_p\cdot |G:C_H(s)|_{p'}$.

The \ir characters of $H$ can be partitioned in classes consisting
of characters of equal degree. This has been done in Nozawa
\cite{No}, who computed the \ir characters of $H$. Similarly, the
elements $g\in H$ can be partitioned in classes consisting of all
elements $g'$ such that $C_H(g')$ is conjugate to $C_H(g)$. Below
we use Nozawa's notation $A_1\ld A_{14}, B_1,\ldots  $ for such
classes. (So each class in question is a union of conjugacy
classes.) In Table 4 below the first column lists the semisimple
conjugacy classes of $H$ with conjugate centralizers $C_H(s)$, and
the second column lists $|C_H(s)|$. In order to extract from
\cite{No} the regular characters we use the above formulas for
their degrees. The third column lists $\rho_s(1)$. For reader's
convenience we also identify $\rho_s$ with notation in \cite{No}
in the fourth column. For instance, $\chi_{11}(s)$ for $s\in A_1$
are characters of the same degree $p^6$ depending on a parameter
$s$. In computations below we do not use $s$ to parameterize the
characters;  instead we write $\chi_{11}(1)=p^6$ to tell that {\it
every} character from the set $\chi_{11}$ takes value $p^6$ at
$1\in H$. (Note that the degrees of  regular characters of
$U(4,q)$ are as in Table 4 with $q$ in place of $p$.)

\medskip
  \begin{center}

 TABLE 4: Degrees of the regular characters of $H=U(4,p)$
    \vspace{10pt}

{\small
   \noindent  \begin{tabular}{|l|l|l|c| }

        \hline
         ~~~~~~~~$~~~~~~~~~~~~~~~~~s$& ~~~~~~$~~~~~~|C_H(s)| $
         &$~~~~~~~~~\rho_s(1)$
        &$\rho_s$ in \cite{No}\\
\hline
     $A_1 $  & $~~~~~~~~~|H|$ & $~~ p^6$ & $\chi_{11}(s) $
      \cr
\hline $A_6 $
&$~p^3(p+1)^2(p^2-1)(p^3+1)$&$p^3(p-1)(p^2+1)$&$\chi_{13}(s)$\\
 \hline
$A_9 $& $~p^2(p+1)^2(p^2-1)^2$& $~p^2(p^2+1)(p^2-p+1)
$&$\chi_{20}(s)$\\
 \hline
$A_{12} $ & $~p(p+1)^3(p^2-1)$&
$~p(p-1)(p^2-p+1)(p^2+1)$&$\chi_{15}(s)$\\ \hline $A_{14} $&
$(p+1)^4$& $(p-1)^2(p^2-p+1)(p^2+1)$&$\chi_{10}(s)$\\ \hline
 $B_1 $ &$~p(p+1)(p^2-1)^2$& $p(p^2+1)(p^3
+1)$&$\chi_{8}(s)$\\ \hline $B_3 $& $(p+1)^2(p^2-1)$& $
(p-1)(p^2+1)(p^3+1)$&$\chi_{6}(s)$\\ \hline
$C_1$&$~p^2(p^2-1)(p^4-1)$&$p^2(p+1)(p^3+1)$&$\chi_{4}(s)$\\
\hline $C_3 $&$(p^2-1)^2$&$(p+1)(p^2+1)(p^3+1)$&$\chi_{2}(s)$\\
\hline $D_1$ &$(p+1)(p^3+1)$& $(p^2-1)(p^4-1) $&$\chi_{9}(s)$\\
\hline $E_1$& $~p^4-1 $& $(p+1)(p^3+1)(p^2-1)$&$\chi_{5}(s)$\\
\hline
    \end{tabular}}
\end{center}

\bigskip

\begin{lemma}\label{u4q}
Let $G=SU(4,p)$ and $\Phi\neq St$ be a PIM.
Then $\dim\Phi>|G|_p$.
\end{lemma}

Proof. (1)
 Let $V$ be  the socle and $\chi$ the character of
$\Phi$. Then  $V=V_\mu|_G$, where $V_\mu$ is a $\GC$-module with
$p$-restricted \hw $\mu$.  Suppose that $c_\Phi=1$; then Lemma
\ref{sd7}  implies that $\mu=(p-1)(\om_1+ \om_3)$.

 We have explained in Section \ref{74} that
$\chi=\tau+\gamma +\sum \nu_i$, where $\tau$ is constituent of
$1_U^G$, $\gamma$ is a regular character and $\nu_i$ are some \ir
characters that are neither regular nor in $1_U^G$. In particular,
  $\nu_i(1)\leq \chi(1)-\tau(1)-\gamma(1)$. We
keep notation of Section \ref{74}.

\med
(2)  By Lemma \ref{de4}, $\Phi=\Psi|_G$, where $\Psi$ is some PIM
for $H$. In particular, $\dim \Psi =p^6=|H|_p$.
Let $\chi$ be the character of $\Psi$, so $\chi(1)=p^6$. 
By Corollary \ref{gg1}, $\chi$ must contain exactly one regular
character $\rho$.

\medskip
(3) If  $(\chi,\gamma)>0$ for a regular character $\gamma$ of $H$
then $\gamma(1)=p(p-1)(p^2+1)(p^2-p+1)$ or
$(p-1)^2(p^2+1)(p^2-p+1)$.

\med
Indeed, set $f=p^6-\tau(1)=p^3(p-1)(p^2+1)$. Then
$\rho(1)\leq f$. One easily checks that for characters $\rho_s$ in
Table 4 $\rho_s(1)>f$ unless $s\in\{A_{6},A_{12},A_{14}\}$.

Observe that $\rho_s(1)=f$ for $s\in A_{6}$. However, $\gamma$
cannot coincide with $\rho_s$ for this $s$, as otherwise
$\tau=\tau'|_G$ for some character $\tau'\in\Irr H$ and
$\chi=\tau'+\rho_s$;  by  inspection in \cite{No}, there are
non-semisimple elements $g\in G\subset H$ such that
$\tau(g)+\rho_s(g)\neq 0$ for $s\in A_{6}$, while $\chi(g)=0$ as
$\chi$ is the character of a projective module. Thus,
$s\in\{A_{12},A_{14} \} $, so $\rho_s\in \{\chi_{15},\chi_{10}\}$,
and (3) follows.

\med

(4) Thus, $\gamma\in\{\chi_{10},\chi_{15} \} $.  Note that
$\chi_{10}(1)<\chi_{15}(1)$. Set

\med
$e_1=f-\chi_{15}(1)=p(p^2+1)(p-1)^2$ and

\med
$e_2=f-\chi_{10}(1)=(p^2+1) (p-1)(2p^2-2p+1)$.

\med
 Then  $\nu_i(1)\leq e_1$ or $e_2$. 
It follows from  Lemma \ref{u44}(3) that $\cup_{s\in T_0^*
}{\mathcal E}_s= \Irr 1_U^G$. Therefore, $\nu_i\in {\mathcal E}_s$
for some semisimple elements $s\in H^*$ for   $s\notin T_0^*$.

We recall some facts of character theory of Chevalley groups.

\med
Observe that  $A_1,A_9,B_1,C_1,C_3$ are the only conjugacy classes
that meet  $T_0$. So $\nu\in {\mathcal E}_s$ and $s\notin
\{A_1,A_9,B_1,C_1,C_3\}$.  It is well known that $ |{\mathcal
E}_s|=1 $ \ii  then $(|C_H(s)|,p)=1$ (this follows for instance
from the formulas for regular and semisimple characters in
${\mathcal E}_s$, see \cite[Ch. 8]{Ca}). So, if $(|C_H(s)|,p)=1$
then the regular character is the only character of ${\mathcal
E}_s$. This happens \ii $s\in \{A_{14}, B_3,C_3,D_1,E_1\}$, see
Table 4. As $\rho$ is the only regular character that is a
constituent of $\chi$, in our case $s\notin \{A_{14},
B_3,C_3,D_1,E_1\}$. We are left with the cases
$s\in\{A_6,A_{12}\}$. Let $S$ denote the subgroup of $C_H(s)$
generated by unipotent elements. By the Deligne-Lusztig theory,
the characters in ${\mathcal E}_s$ are of degree  $d\cdot
|H:C_H(s)|_{p'}$, where $d$ is the degree of a unipotent character
of $S$. If $s\in A_{12}$ then $S\cong SL(2,p)$, and if $s\in A_6$
then $S\cong SU(3,p)$, see \cite{No}. The degrees of unipotent
characters of these groups are well known to be $1,p$ and
$1,p(p-1),p^3$, respectively. Therefore, non-regular characters in
${\mathcal E}_s$ are of degrees $\chi_{19}(1)=(p-1)(p^2+1)$ and
$\chi_{17}(1)=p(p-1)^2(p^2+1)$ for $s\in  A_6$, and of degree
$\chi_{16}(1)=(p-1)(p^2-p+1)(p^2+1)$ for $s\in A_{12}$. Note that
$\chi_{16}(1)> \chi_{17}(1)=e_1>\chi_{19}(1)$.

\med
For further use, we write down some character values extracted
from \cite{No}. Let  $g\in A_{10}$, $h\in A_{11}$ with the same
semisimple parts, and $\chi_i$ are as in the last column of Table
4.

\med
$~~~~~~~~~~~~~(*) \,\,\,\,\,\,\,\,\,\,\,\,\,\,\, \chi_i(g)+(p-1)\chi_i(h)=
\begin{cases}0& for ~~i=10~~and ~~i=17\\
\pm p&for ~~i=13\\
\pm 1&for ~~i=16~~and ~~i=19\\
 \end{cases}
$

\med
Note that the absolute value of $\chi_i(g)+(p-1)\chi_i(h)$ is
independent from the choice of an individual  character in   the
set $\chi_i$. This is the reason to consider
$\chi_i(g)+(p-1)\chi_i(h)$ instead of computing $\chi$ at $g$ or
$h$  in some formulas below.

\med
 Case 1. $\gamma (1)=\chi_{15}(1)$.

\med
The \rep  $\tau$ above is denoted by $\chi_{13}$ in \cite{No}.
Note that $\chi_{13}(1)+\chi_{15}(1)+\chi_{17}(1)=p^6$. However,
$\chi$ is not the sum of characters from these sets as
$\chi_{13}+\chi_{15}+\chi_{17}$ does not vanish at elements of
some class in $A_{11}$.


\med
It follows that the only possibility is
$\chi=\chi_{13}+\chi_{15}+p(p-1) \chi_{19}$. (This is not an
actual formula, it only  tells that $\chi$ is the sum of $p(p-1)$
characters from the set $\chi_{19}$ and one character from each
set $\chi_{13}$ and $\chi_{15}$.) Inspection of \cite{No} shows
that this is false. (One can compute the values of these
characters at a regular unipotent element of $H$; it takes zero
values for every character from $\chi_{13}$ and $\chi_{15}$ and
value $-1$ for every character from $\chi_{19}$. So we have a
contradiction.)

\med
 Case 2. $\gamma(1)=\chi_{10}(1)$.

\med
Then $e_2=(p^2+1) (p-1)(2p^2-2p+1)=x\cdot \chi_{16}(1)+y\cdot \chi_{17}(1)+z\cdot \chi_{19}(1)$

\med $=x(p-1)(p^2-p+1)(p^2+1)+yp(p-1)^2(p^2+1)+z(p-1)(p^2+1)$,

\med
where $x,y,z\geq 0$ are integers. Dividing by $(p^2+1)(p-1)$, we get:
\med

$2p^2-2p+1=x(p^2-p+1)+yp(p-1)+z$, or $p(p-1)(2-x-y)=z+x-1$.

\med
An obvious possibility is  $x=y=0$, $z=2p^2-2p+1$. Suppose $x+y\geq 1$.
If $x+z=1$
then either $z=0, x=y=1$ or $z=1,x=0,y=2$. 
\med

Suppose $x+z\neq 1$. Then $0\neq x+z-1\equiv 0\mod p(p-1)$, so $x+z-1=kp(p-1)$ for some integer $k>0$.
This implies $2=x+y+k$, so $x+y\leq 1$, and hence $x+y=1$, $k=1$. So either
  $x=1,y=0$, $z=p(p-1)$ or $x=0,y=1,z=p(p-1)+1$.

Altogether we have five solutions.

\med
(1) $e_2=(2p^2-2p+1)\chi_{19}(1)$, and $\chi=\chi_{13}+\chi_{10}+(2p^2-2p+1)\chi_{19}$;

\med
(2) $e_2=\chi_{16}(1)+\chi_{17}(1)$ and $\chi=\chi_{13}+\chi_{10}+\chi_{16}+\chi_{17}$;

\med
(3) $e_2=2\chi_{17}(1)+\chi_{19}(1)$ and $\chi=\chi_{13}+\chi_{10}+2\chi_{17}+\chi_{19}$;

\med
(4)  $e_2=\chi_{16}(1)+(p^2-p)\chi_{19}(1)$ and $\chi=\chi_{13}+\chi_{10}+\chi_{16}+(p^2-p)\chi_{19}$;

\med
(5) $e_2=\chi_{17}(1)+(p^2-p+1)\chi_{19}(1)$ and $\chi=\chi_{13}+\chi_{10}+\chi_{17}+(p^2-p+1)\chi_{19}$.

\med
As above, here $(2p^2-2p+1)\chi_{19}$ means the sum of $(2p^2-2p+1)$ characters from the set $\chi_{19}$, and similarly for other cases. 

\med
 Let  $g\in A_{10}$, $h\in A_{11}$ with the same semisimple
parts. Next, we compute $\chi(g)-(p-1)\chi(h)$ for $\chi$ in the
cases (2) and (3) above. As both $g,h$ are not semisimple, this
equals 0. However, computing this for the right hand side in these
formulas gives us a non-zero value. This will rule out these two
cases.

\med
Case 2: Using (*), we have $0=\chi(g)+(p-1)(\chi(h)=\chi_{10}(g)+\chi_{13}(g)+\chi_{16}(g)+\chi_{17}(g)
+(p-1)(\chi_{10}(h)+\chi_{13}(h)+\chi_{16}(h)+\chi_{17}(h)) =
 \pm p\pm 1\neq 0$. 
This is a contradiction. 

\med
Case 3: Similarly, by (*), we have  $0=\chi(g)+(p-1)\chi(h)=\chi_{10}(g)+\chi_{13}(g)+2\chi_{17}(g)+\chi_{19}(g)+
(p-1)(\chi_{10}(h)+\chi_{13}(h)+2\chi_{17}(h)+\chi_{19}(h))=\pm p\pm 1\neq 0$.

\med
In order to deal with cases (1), (4), (5), we compute the character values at the regular unipotent element $u$ (from class $A_5$ of \cite{No}).
We have $\chi_{10}(u)=1$, $\chi_{13}(u)=0$, $\chi_{16}(u)=-1$,
$\chi_{17}(u)=0$, $\chi_{19}(u)=1$. The values do not depend on the choice of an individual character in every class $\chi_i$ ($i\in \{10,13,16,17,19\}$.)

\med
Case (1). $\chi(u)=\chi_{10}(u)+\chi_{13}(u)+(2p^2-2p+1)\chi_{19}(u)=
1-(2p^2-2p+1)\neq 0$, a contradiction.

\med
Case (4). $\chi(u_5)=\chi_{10}(u_5)+\chi_{13}(u_5)+\chi_{16}(u_5)+(p^2-p)\chi_{19}(u_5)=
1-1-(p^2-p)\neq 0$, a contradiction.

\med
Case (5).
$\chi(u_5)=\chi_{10}(u_5)+\chi_{13}(u_5)+\chi_{17}(u_5)+(2p^2-2p+1)\chi_{19}(u_5)=
1-(2p^2-2p+1)\neq 0$, a contradiction. \hfill $\Box$

\medskip
Remark. Some characters $\chi$ in cases (2) and (3)  vanish at
non-identity $p$-elements. (In (3) $2\chi_{17}$ may be the sum of
two distinct characters of degree $p(p-1)^2(p^2+1)$.)

\subsection{The groups $G={}^3D_4(p)$}

We follows the strategy described in Section \ref{74}. We can
assume $p>2$ as the decomposition matrix for ${}^3D_4(2)$ is
available in the GAP library, and one can read off from there that
the minimum valu e for $c_\Phi$ equals 15.

\begin{lemma}\label{d4x}
Let $G={}^3D_4(p)$, and $\Phi\neq St$
be a PIM. Then $c_\Phi\geq 2$.
\end{lemma}

Proof. Let $V$ be the socle and $\chi$ the character of $\Phi$.
Then
$V=V_\mu$, where $\mu=(p-1)(\om_1+\om_3+\om_4)$. 
We keep notation of Section \ref{74}.

Arguing by contradiction,  we assume that $c_\Phi=1$, and then we
denote by $\tau$ and $\gamma$ \ir characters of $G$ occurring as
common constituents of $\chi$ with
 $1_U^G$ and $\Gamma$, respectively.

\med

We first show that $\tau$ is a unipotent character of degree
$p^7(p^4-p^2+1)$. Let $P$ be a parabolic subgroup of $G$
corresponding to the nodes  $J=\{1,3,4\}$ at the Dynkin diagram of
$G$. Let $L$ be a Levi subgroup of $P$ and $L'$ the subgroup of
$L$ generated by unipotent elements. Then $L'\cong SL(2,p^3)$. Let
$S=\soc (V_\mu|_P)|_L$ and $\rho$ the character of $S$. By Lemma
\ref{nf1}, $(\chi, \rho^{\# G})>1$. As $\rho^G$ is a part of
$1_U^G$, it follows that $(\chi, \rho^G)=1$. In addition, $\tau$
is a constituent of $1_B^G$. The degrees of \ir constituents of
$1_B^G$ are given in \cite[Proposition 7.22]{CKS}. The order of
$G$ is $p^{12}(p^6-1) (p^2-1)(p^8+p^4+1)$, so
 $\rho^{\# G}(1)=p^3|G:P_J|=p^3(p+1)(p^8+p^4+1)$. From this one easily
obtains the \f lemma (which true for $q$ in place of $p$):

\begin{lemma}\label{d41}
In the above notation,
 $\rho^{\# G}=St+\rho'_1+\rho'_2+\rho_2$, where
$ \rho_2'(1)=p^3(p^3+1)^2/2$, $ \rho_2(1)=p^3(p+1)^2(p^4-p^2+1)/2$
and $ \rho_1'(1)=p^7(p^4-p^2+1)$.
\end{lemma}

\begin{lemma}\label{d42}
 Let $\Phi$ be a PIM with $c_{\Phi}=1$.
Then $\tau(1)=p^7(p^4-p^2+1)$.
\end{lemma}

Proof. By \cite[Proposition 7.22]{CKS} or \cite[p.115]{CKS}, the
characters $\rho_2,\rho'_2$ occur in $1_B^G$ with multiplicity
greater than 1; therefore, none of them is a constituent of
$\chi$, and the claim follows.\hfill $\Box$

\med
As  $G$ coincides with its dual group $G^*$, we identify maximal
tori in $G$ and $G^*$. Following \cite{DeM} we denote by $s_i$,
$i=1\ld 15$, the union of the semisimple classes of $G$ with the
same centralizer (that is, $C_G(x)$ is conjugate with $ C_G(y)$
for $x,y\in s_i$). In the character table of $G$ in \cite{DM} a
class with representative $s\in s_i$ meets $ T_0$ \ii $i\leq 8$.
(The set $s_8$ consists of regular elements.)

\begin{lemma}\label{uc4} The set $\cup _{i\leq 8}{\mathcal E}_{s_i}\setminus \Irr 1_U^G$
consists of two unipotent cuspidal characters.\end{lemma}

  Proof.  The unipotent characters of $G$ have been determined by Spaltenstein
  \cite{Spa}. All but two of them are constituents of $1_B^G$.

   The characters
$\chi_{3,*}$, $\chi_{5,*}$, $\chi_{6,*}$, $\chi_{7,*}$ and
$\chi_8$ from $\cup _{i\leq 8}{\mathcal E}_{s_i} $ are either
regular or semisimple. So they are in $ 1^G_U $ by Lemma
\ref{dd4}(2). This  leaves with $s_i$ for $i=2,4$ (as $s_1=1$). In
these cases $C_G(s_2)\cong (SL(2,p^3)\circ SL (2,p))\cdot  T_0$
and $C_G(s_4)\cong  SL (3,p)\cdot  T_0$, where $T_0$ is a split
maximal torus in $G$, see \cite[Proposition 2.2]{DeM}. For these
groups the number of unipotent characters are known to be 4 and 3,
respectively. So $|{\mathcal E}_{s_i} |=4$, resp., 3, for $i=2,4$,
see \cite[Theorem 13.23]{DM}. Let $\beta_i\in\Irr T_0$ be the
character corresponding to $s_i$. Then the Harish-Chandra series
$\beta_{i,B}^G$ is contained in $ {\mathcal E}_{s_i} $ by Lemma
\ref{dd4}(1). Set $W_i=C_{W_0}(\beta_i)$. Then $W_2\cong
\ZZ_2\times \ZZ_2$ and $W_4\cong S_3$ \cite[Lemma 3.4]{DeM}.
Moreover, the centralizing algebra of $\beta_{i,B}^G$  is
isomorphic to the group algebra of $W_i$ (\cite[Exercise in \S 14
prior Lemma 86]{St}). Therefore, $\beta_{1,B}^G$ consists of 4
distinct \ir constituents, whereas $\beta_{2,B}^G$ has three
 distinct \ir constituents (two of them occurs with \mult 1 and one
 constituent occurs with  \mult 2). In both the cases $|{\mathcal
E}_{s_i} |$ coincides with the number of distinct \ir constituents
 in $\beta_{i,B}^G$, and the claim follows. \hfill $\Box$

\begin{lemma}\label{d49}
$\gamma(1)\in\{d_1,d_2,d_4\}$, where
$d_2=(p-1)^2(p^3+1)^2(p^4-p^2+1)$, $d_4=(p-1)(p^3-1)(p^8+p^4+1)$.
\end{lemma}

Proof. We first write down the degrees of the regular characters
of $H$ in Table 5. Note that one can easily detect the regular
characters in the character table of $G$ in \cite{DeM}. Indeed,
the characters in \cite{DeM} are partitioned in Lusztig series
${\mathcal E}_{s}$ with $s\in
s_i$, $i=1\ld 15$. 
As mentioned in the proof of Lemma \ref{gg4}, 
every ${\mathcal E}_s$ has a single regular character. Then its
degree equals $|G:C_G(s)|_{p'}\cdot |C_G(s)|_p$, while other
characters in ${\mathcal E}_s$ are of degree $|G:C_G(s)|_{p'}\cdot
e$, where $e$ is the degree of a unipotent character in $C_G(s)$.
(So the $p$-power part of the character degrees in ${\mathcal
E}_s$ is maximal for the regular character.)

As is explained in Section \ref{74}, $\gamma (1)\leq p^{12}-
\tau(1)=p^{12}-p^{11}+p^{9}-p^{7}$. One observes  that for $s\in
\cup_{i>8}s_i$ only the characters  $\chi_{12}(s)$ and
$\chi_{15}(s)$ satisfy this inequality. So exactly one of these
characters is a constituent of $\chi$, as claimed.

\med

Let $d_1,d_2$ be the degrees of these $\chi_{12}(s)$,
$\chi_{15}(s)$, respectively, and set
$f_i=p^{12}-p^{11}+p^{9}-p^{7}-d_i$ ($i=1,2$). Then

\med
$f_1=p^{11}+p^9+4p^8+p^7  +2p^6+ 2p^5+4p^4+2p-1$,
\med

$f_2=2p^{9}-2p^8+p^5-2p^4+p^3+p-1=(p-1)(2p^8+p^4-p^3+1).$

\med
Thus, $f_i$ $(i=1,2)$ is a sum of the degrees of \ir non-regular
characters that do not belong to $1_U^G $. Let $\lam$  be one of
these characters. Note that $G$ has two cuspidal unipotent
characters of degrees $e_1=p^3(p^3-1)^2/2$ and
$e_2=p^3(p-1)^2(p^4-p^2+1)/4$. They do not belong to  $1_U^G $
whereas the other 6 unipotent characters belong to  $1_U^G $.

 Inspection of the
character table of $G$ in \cite{DeM} shows that non-regular
characters that are  in ${\mathcal E}_{s_i}$ with $i>8$ are the
characters $\chi_{9,1}$, $\chi_{9,qs'}$ and $\chi_{10,1}$ in
\cite{DeM}, of degrees $e_3=(p^3-1)(p^2+p+1)(p^4-p^2+1)$,
$e_4=p(p^3-1)^2(p^4-p^2+1)$ and $e_5=(p^3-1)(p^8+p^4+1)$,
respectively. One observes that $(p^2-p+1)(p^2+p+1)=p^4-p^2+1$ and
$(p^2+p+1)(p-1)=p^3-1$. As $p^2+p+1$ is odd, it follows that
$p^2+p+1$ divides $e_i$ for every $i=1,2,3,4,5$, and hence
$p^2+p+1$ divides $f_j$ for $j=1,2$. However, $f_1({\rm mod }\,
p^2+p+1)=-5$ and $f_2({\rm mod }\, p^2+p+1)=2p-11$. This is a
contradiction, which completes the proof. \hfill $\Box$

\bigskip
  \begin{center}

 TABLE 5: Degrees of the regular characters of $G={}^3D_4(p)$, $p>2$
    \vspace{10pt}

{\small
   \noindent  \begin{tabular}{|l|l|c| }

        \hline
         ~~~~$s$
         &$ \, \, \, \, \, \, \, \, \, \, \, \, \, \, \, \, \,
          \, \, \, \, \, \, \, \, \, \, \, \, \, \, \,\rho_s(1)$
        &$\rho_s$ in \cite{DM}\\
\hline
     $s_1 $  &   $~~ p^{12}$ & $St $
      \cr
\hline $s_2 $
&$~p^4(p^8+p^4+1)
$&$\chi_{2,St,St'}(s)$\\

 \hline
$s_3 $
 & $~p^3(p+1)(p^8+p^4+1)
$&$\chi_{3,St}(s)$\\

 \hline
$s_{4} $ &   $p^3(p^3+1)(p^2-p+1)(p^4-p^2+1)
$&$\chi_{4,St}(s)$\\

\hline $s_{5} $&   $p(p^3+1)(p^8+p^4+1)
$&$\chi_{5,St}(s)$\\

\hline
 $s_6 $ & $(p+1)(p^3+1)(p^8+p^4+1)
$&$\chi_{6}(s)$\\

\hline $s_7 $& $ p^3(p-1)(p^8+p^4+1)$&$\chi_{7,St}(s)$\\

\hline $s_8$& $(p-1)(p^3+1)(p^8+p^4+1)$&$\chi_{8}(s)$\\

\hline $s_9 $& $p^3(p^3-1)(p^2+p+1)(p^4-p^2+1)$&$\chi_{9,
St}(s)$\\

\hline $s_{10}$& $p(p^3-1)(p^8+p^4+1) $&$\chi_{10,St}(s)$\\

\hline $s_{11}$& $(p+1)(p^3-1)(p^8+p^4+1)$&$\chi_{11}(s)$\\

\hline $s_{12}$& $(p-1)^2(p^3+1)^2(p^4-p^2+1)$&$\chi_{12}(s)$\\

\hline $s_{13}$& $(p+1)^2(p^3-1)^2(p^4-p^2+1)$&$\chi_{13}(s)$\\

\hline $s_{14}$&   $(p^6-1)^2$&$\chi_{14}(s)$\\

\hline $s_{15}$& $(p-1)(p^3-1)(p^8+p^4+1)$&$\chi_{15}(s)$\\

\hline
    \end{tabular}}
\end{center}

\section{The minimal PIM's}

In this section we complete the proof of Theorem \ref{wt1}. To
settle the base of induction, we refer to certain known results
for groups of small rank. We first write down some data available
in the GAP library:

\begin{lemma}\label{dn1}
Let $G$ be one of the groups below and let $\Phi$ be a PIM for $G$
other than the Steinberg PIM.

\med

$(1)$ If $G=Sp(4,2)$ or $Sp(4,3)$ then $c_{\Phi}\geq  3$.

\med
$(2)$ If $G=Sp(4,3)$ then $c_{\Phi_1}=2$, and $c_{\Phi} \geq 3$
for  $\Phi\neq \Phi_1$.

\med
$(3)$  If $G=SU(4,2)$ then ${\rm min}\, c_\Phi\geq 4$.

\med
$(4)$  If $G=SU(5,2)$ then ${\rm min}\, c_\Phi\geq 5$.

\med
$(5)$  If $G={}^3D_4(2)$ then $ c_\Phi\geq 15$.

\med
$(6)$  If $G=G_2(2)$ then ${\rm min}\, c_\Phi=5$.

\med
$(7)$  If $G={}^2F_4(2)$ then $ c_\Phi \geq 14 $.

\end{lemma}

\begin{lemma}\label{sm6} Let $G$ be one of the groups below and let $\Phi\neq St$ be a PIM
for $G$.

\med
$(1)$ Let $G=SL(3,p)$, $p>2$. Then $ c_\Phi \geq 2 $. 


\med
$(2)$ Let $G=Sp(4,p)$, $p>3$. Then $ c_\Phi \geq 3 $.

\med
$(3)$ Let $G=G_2(p)$, $p>2$. Then $ c_\Phi \geq 6 $.
\end{lemma}

Proof. The statements (1),  (2), (3) follows from the results in
\cite{Hu81},  \cite{Hu80}, \cite[Ch. 18]{Hub}, respectively.
\hfill $\Box$

\med

{\bf Proof of Theorem} \ref{wt1}. Let $r$ be the BN-pair rank of
$G$. If $r=1$  then the result is contained in Proposition
\ref{ra1}. In general, suppose the contrary, and let $G$ be a
counter example, so $r>1$. Let $\Phi\neq St$ be a PIM with
$c_\Phi=1$, and let $V$ be the socle of $\Phi$. Then $V$ is \ir
and hence $V=V_\mu|_G$ for some \ir module $V_\mu$ for the
respective algebraic group  $\GC$. (Here $\mu$ is the highest
weight of $V_\mu$.) Let $P$ be a proper parabolic subgroup of $G$,
which is not a Borel subgroup of $G$. Then $P=P_J$ for some
non-empty set $J$, see Section 5. Let $L_J$ be a Levi subgroup of
$P$,  let $G_J$ be as in Section 5 and $V_J=C_V(O_p(P))$. Then
$G_J$ is a Chevalley group of rank $r_J<r$, and $V_J$ is  \ir both
as an $FL_J$- and an $FG_J$-module (Lemma \ref{le1}). Let $\Psi$
be the PIM for $L_J$ with socle $V_J|_{L_J}$. By Lemma \ref{mp1},
$c_\Psi=1$, so $\dim \Psi=|L_J|_p$. As $|L_J|_p=|G_J|_p$, the
restriction $\Psi|_{G_J}$ is a PIM for $G_J$ of dimension
$|G_J|_p$ (Lemma \ref{om1}(2)). This is a contradiction unless
$V_J|_{G_J}$ is isomorphic to the Steinberg module for $G_J$, or
else $G_J\in \{A_1(p), A_2(2),{}^2A_2(2)\}$ and $V_J|_{G_J}$ is
the trivial $FG_J$-module. As this is true for every parabolic
subgroup of $G$, which is not a Borel subgroup, $V$ satisfies the
assumption of Lemma \ref{sd7}. Therefore, $G$ belongs to the list
(1) - (6) of Lemma \ref{sd7}, or else $G={}^2F_4(\sqrt{2})$.
However, for these groups  Theorem \ref{wt1} is true by Lemmas
\ref{u4q}, \ref{d4x}, \ref{dn1} and \ref{sm6}, which is a
contradiction.

\medskip
Acknowledgement. The author is indebted Ch. Curtis, J. Humphreys,
G. Malle and M. Pelegrini for discussion and remarks.

\end{document}